\documentclass[nonblindrev]{informs3} 

\OneAndAHalfSpacedXI 

\usepackage{endnotes}
\let\footnote=\endnote
\let\enotesize=\normalsize
\def\notesname{Endnotes}%
\def\makeenmark{$^{\theenmark}$}
\def\enoteformat{\rightskip0pt\leftskip0pt\parindent=1.75em
  \leavevmode\llap{\theenmark.\enskip}}

\usepackage{natbib}
 \bibpunct[, ]{(}{)}{,}{a}{}{,}%
 \def\bibfont{\small}%
\TheoremsNumberedThrough     
\ECRepeatTheorems

\EquationsNumberedThrough    

\usepackage[utf8]{inputenc}
\usepackage{multirow}
\usepackage[font=scriptsize]{caption}
\usepackage{subcaption}
\usepackage{algorithm}
\usepackage{algorithmic}
\usepackage{xr-hyper}
\usepackage[hidelinks]{hyperref}
\hypersetup{
     colorlinks = true,
     citecolor  = blue,
     linkcolor = blue,
     urlcolor = blue
}
\usepackage{url}
\usepackage{color}
\usepackage{siunitx}
\usepackage{bbm}

\newcommand{\mynameref}[1]{\textcolor{blue}{\textit{\nameref{#1}}}}
\DeclareMathOperator*{\lexmin}{lex-min}
\usepackage{enumitem}
\usepackage{booktabs}
\usepackage{soul}
\usepackage{pifont}

\newtheorem{properties}{Property}
\newtheorem{thm}{Theorem}

\begin{document}


\RUNAUTHOR{Guan, Basciftci, and Van Hentenryck}

\RUNTITLE{Path-Based Formulations for ODMTS-DA Problem}

\TITLE{Path-Based Formulations for the Design of On-demand Multimodal Transit Systems \\ with Adoption Awareness}

\ARTICLEAUTHORS{%
\AUTHOR{Hongzhao Guan}
\AFF{H. Milton Stewart School of Industrial and Systems Engineering, Georgia Institute of Technology, \EMAIL{hguan7@gatech.edu}} 
\AUTHOR{Beste Basciftci}
\AFF{Department of Business Analytics, Tippie College of Business, University of Iowa, \EMAIL{beste-basciftci@uiowa.edu}}
\AUTHOR{Pascal Van Hentenryck}
\AFF{H. Milton Stewart School of Industrial and Systems Engineering, Georgia Institute of Technology, \EMAIL{pascal.vanhentenryck@isye.gatech.edu}}
} 

\ABSTRACT{%

This paper reconsiders the On-Demand Multimodal Transit Systems (ODMTS) Design with Adoptions problem
(ODMTS-DA) to capture the latent demand in on-demand multimodal
transit systems. The ODMTS-DA is a bilevel optimization problem, for
which \citet{Basciftci2021} proposed an exact combinatorial Benders
decomposition. Unfortunately, their proposed algorithm only finds
high-quality solutions for medium-sized cities and is not practical
for large metropolitan areas.
The main contribution of this paper is to propose a new path-based
optimization model, called {\sc P-Path}, to address these computational
difficulties. The key idea underlying {\sc P-Path} is to enumerate two
specific sets of paths which capture the essence of the choice model associated with the adoption behavior of riders.
With the help of these path sets, the ODMTS-DA can be formulated as a
single-level mixed-integer programming  model. In addition, the paper presents preprocessing
techniques that can reduce the size of the model significantly.  {\sc
  P-Path} is evaluated on two comprehensive case studies: the
mid-size transit system of the Ann Arbor -- Ypsilanti region in
Michigan (which was studied by \citet{Basciftci2021}) and the large-scale transit system for the city of
Atlanta.  The experimental results show that {\sc P-Path} solves the
Michigan ODMTS-DA instances in a few minutes, bringing more than two orders of magnitude improvements compared to the existing
approach. For Atlanta, the results show that {\sc P-Path} can solve
large-scale ODMTS-DA instances (about 17 millions variables and 37
millions constraints) optimally in a few hours or in a few
days. These results show the tremendous computational benefits of {\sc
  P-Path} which provides a scalable approach to the design of on-demand
multimodal transit systems with latent demand.
}%



\maketitle

%



\section{Introduction}

\begin{figure}[htb]
    \centering
    \includegraphics[width=0.8\textwidth]{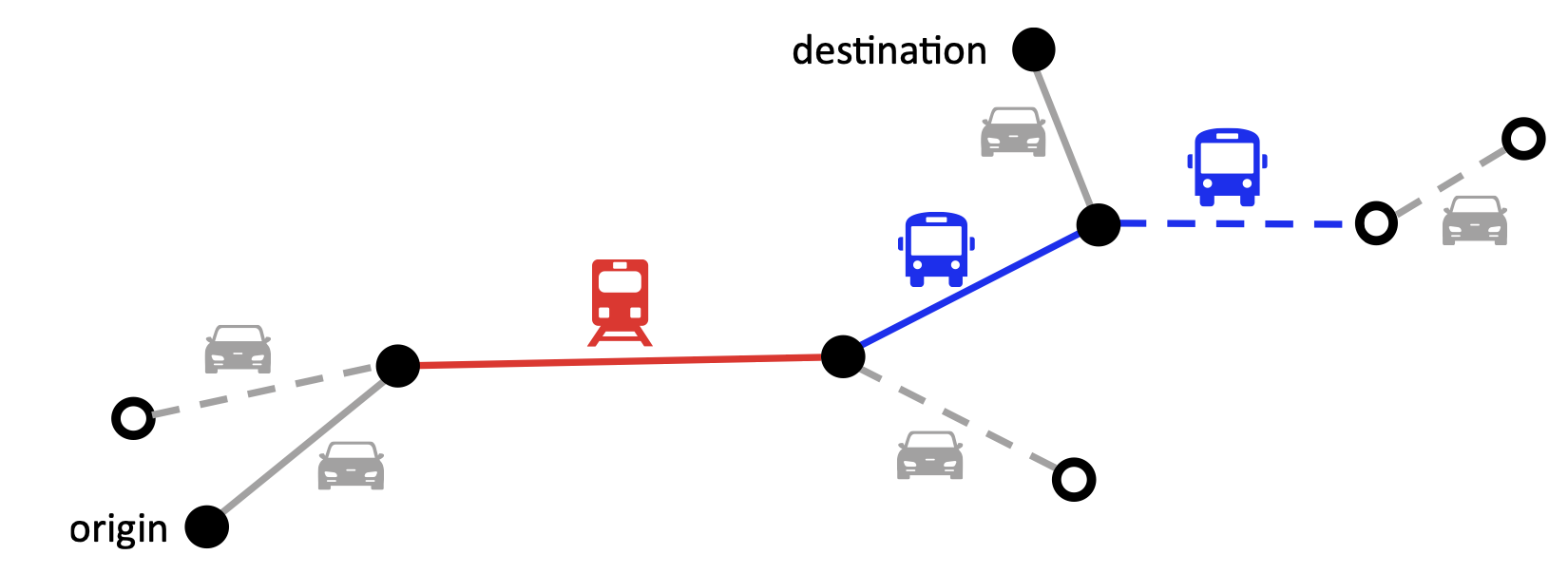}
    \caption{An example of ODMTS with a sample Rider Path from their origin to destination (solid line)}
    \label{fig:ODMTS}
\end{figure}

On-Demand Multimodal Transit Systems (ODMTSs) tightly integrate
on-demand dynamic shuttles with fixed transit services such as rail
and high-frequency buses. 
As illustrated in Figure~\ref{fig:ODMTS}, ODMTSs are operated around a number of hubs
that are transit stations for high-frequency buses and urban
rails. On-demand shuttles in ODMTSs are primarily utilized as the
feeders to/from the hubs and represent an effective solution to the
first-mile and last-mile problem faced by the vast majority of transit
agencies. A realistic ODMTS pilot conducted in Atlanta, Georgia, USA, has demonstrated the efficacy of ODMTS in providing efficient services and economical sustainability. ODMTSs particularly benefit local communities engaging in short-distance trips or seeking connections to high-frequency routes, such as the rail system in Atlanta \citep{van2023marta}.
Furthermore, several simulation studies have shown that ODMTSs can provide
significant cost and convenience benefits under varied settings
\citep{maheo2019benders, dalmeijer2020transfer,
  auad2021resiliency}. However, these studies focus on designing an
ODMTS for the existing transit users and neglect additional riders who
could potentially adopt the system, given its higher convenience.  


To address this gap by integrating rider behavior, \citet{Basciftci2021} propose the ODMTS Design with
Adoption (ODMTS-DA) problem that captures the latent demand in
ODMTS. The ODMTS-DA is a bilevel
optimization problem, where the transit agency suggests riders paths from their origins to destinations by utilizing the transit network design and on-demand shuttles. After observing this path, the riders have the choice to adopt or reject the ODMTS. 
To solve this bilevel problem, an exact combinatorial
Benders decomposition was proposed by taking into account the binary nature of the rider choices. Nevertheless, their proposed algorithm only
finds high-quality solutions for medium-sized cities and is not
practical for large metropolitan cities.

This paper reconsiders the ODMTS-DA and proposes a new path-based
optimization model, called {\sc P-Path}, to address these computational
difficulties. The main contributions of this paper can be summarized as follows:
\begin{enumerate}
\item It proposes {\sc P-Path}, a new path-based optimization model
  that replaces the bilevel subproblems by computing a number of
  specific path sets, that capture the essence of the subproblems and the mode choice models. 
  By leveraging these path sets, the bilevel model can then be
  reformulated as a single-level mixed-integer programming (MIP) model.
\item It presents a number of preprocessing techniques to reduce the
  sizes of the path sets and the number of variables and constraints in the MIP model, enhancing computational efficiency.
\item It demonstrates, on two case studies, the computational benefits of {\sc
  P-Path},  which provides a scalable approach to
  the design of on-demand multimodal transit systems with latent demand. The mid-size transit system of the Ann Arbor -- Ypsilanti region in
  Michigan (which was studied by \citet{Basciftci2021}) and the large-scale transit system for the city of Atlanta. The computational experiments show that {\sc P-Path} solves different Michigan ODMTS-DA instances in a few minutes, demonstrating more than two orders of magnitude improvements in comparison to the existing approach.   The results over various Atlanta ODMTS-DA instances (comprising over 50,000 trips) show that {\sc P-Path} can obtain optimal solutions for large-scale instances, with up to 17 millions of variables and 37 millions of constraints derived from a graph with more than 2,400 nodes, within a few hours/days.
\end{enumerate}

\noindent To ease understanding, {\sc P-Path} is presented in two
steps. In a first step, another path formulation, called {\sc C-Path},
is presented. {\sc C-Path} replaces the follower subproblems in the original bilevel ODMTS-DA formulation by
enumerating all the paths: it also leads to a single-level MIP model
but unfortunately this formulation is not practical given the massive
number of paths even in mid-size instances. In a second step, the
presentation turns to the key innovations in {\sc P-Path}. First, {\sc
  P-Path} only enumerates two specific sets of paths: paths that are
adopted by riders given their choice models and paths that are
profitable for the transit agency but rejected by the riders. Second,
{\sc P-Path} reasons about paths implicitly in the MIP formulation by
using arc variables. The two path sets enumerated by {\sc P-Path} are
small in size, which makes the path enumeration and the MIP model
tractable even for large-scale instances.

The remaining part of this paper proceeds as follows: Section
\ref{sect:literature} presents the relevant literature. Section
\ref{sect:ProblemStatement} introduces the ODMTS-DA formulation.
Section~\ref{sect:paths_method} introduces the proposed path-based
formulations: {\sc C-Path} and {\sc P-Path}. Section
\ref{sect:computational_enhancements} presents preprocessing
techniques to reduce the size of the path sets and the MIP
formulation.  Sections~\ref{sect:ypsi_case_study} and
\ref{sect:atl_case_study} consider real case studies on medium-sized
and large-scale instances, respectively, to demonstrate the
performance of {\sc P-Path}. Section~\ref{sect:conclusion} concludes
the paper.


\section{Literature Review}
\label{sect:literature} 

Transit network design is an important problem in urban planning and
transportation that produces a connected network for cities
\citep{FARAHANI2013281}. 
Examples of transit network design on traditional transit modes, such as buses and rails, can be found in various works such as \cite{borndorfer2007column, fan2018optimal}, and \cite{almasi2021urban}. The emergence of Mobility-as-a-Service (MaaS) over the past decade has led to a substantial transformation of the transit landscape \citep{shaheen2016mobility}. Consequently, the transit network design problem has garnered increased attention with the integration of MaaS services \citep{stiglic2018enhancing, pinto2020joint, liu2021mobility, najmi2023multimodal}.

ODMTSs are transit systems that integrate on-demand services and conventional transit modes. A key distinction between ODMTSs and the majority of MaaS systems lies in the fact that ODMTSs are entirely operated by transit agencies, as opposed to on-demand services provided by Transit Network Companies (TNCs) \citep{van2023marta}. In this regards, the core concepts underlying
the emerging ODMTSs can be traced back to the multi-modal hub location
and hub network design problems \citep{alumur2012multimodal}. Those
were extended by ODMTSs to employ the on-demand shuttles for addressing
the first-and-last-mile problem, and establishing hub-to-hub bus
networks where shuttles are used as feeders
\citep{maheo2019benders}. Continuing from this line of work,
\citet{dalmeijer2020transfer} incorporated bus arc frequencies,
passenger transfer limit, and backbone transit lines such as existing
rail services. These ODMTS studies consider a fixed transit demand
and aim to discover a network design that minimizes a weighted
combination of the passenger convenience and operating cost from the
transit agency perspective. However, neither of these studies take into
account the latent demand when solving the design problem.

Like the first two ODMTS studies, many methods for transit network
design problems are based on the assumption that the passenger demand
is known and fixed \citep{schobel2012line}. Network design problem
with fixed demand can also be generalized to other fields such as
supply chain, energy distribution, and telecommunications. For
instance, \citet{gudapati2022network} present a path-based approach to
solve large-scale network design problem with service requirements. On
the other hand, another group of studies aims at integrating the
transit network design problem together with latent demand, i.e.,
potential customers who might switch to the proposed transit networks
once they are built. For example, instead of solving an optimization
problem with a fixed demand, \citet{klier2008line} modeled the impact
of travel time on transit demand using a linear demand function of the
expected travel time.

A series of studies presented by Canca et al. comprehensively
discussed network designs and latent demand. Their first study
presented a model that concurrently determines the network design,
line planning, capacity, fleet investment, and passengers mode and
route choices \citep{canca2016general}. In the two follow-up studies,
the authors proposed an adaptive large neighborhood search
meta-heuristic for the previous model and a bilevel meta-heuristic
which is designed to solved a revised model with large-scale data,
respectively \citep{canca2017adaptive, canca2019integrated}.

Among the studies that consider latent ridership, three systematic
studies have investigated ODMTS with rider adoption awareness
\citep{Basciftci2020, Basciftci2021, guan2022}. The first two studies
indicate that bilevel optimization frameworks need to be utilized to
model the problem to ensure a fair design for the transit system
\citep{Basciftci2020, Basciftci2021}. In general, a bilevel
formulation is a mathematical program in which several variables are
constrained by the solution of another optimization problem
\citep{kleinert2021survey}. Hence, a bilevel formulation is useful
when modeling hierarchical decision process and has been widely
utilized in the field of transportation. To model adoptions, in the
first study, a personalized choice model which associates adoption
choices with the time and cost of trips in the ODMTS is incorporated
into the bilevel optimization framework \citep{Basciftci2020}. In the
second study, the choice model was then redesigned such that the rider
adoption choices solely depend on the trip duration under a fixed
pricing strategy of the transit agency \citep{Basciftci2021}. These
two studies propose formulations that are called the ODMTS Design with
Adoptions (ODMTS-DA) problem. The authors of these two studies also
provide significant insights into the design of exact algorithms that
decompose and solve the problem under different choice model
assumptions. However, the complexity and combinatorial nature of the
exact algorithms led to computational difficulties: the exact
algorithms can only provide high-quality solutions to medium-sized
ODMTS-DA models.  In order to address this difficulty, another study
further investigated the properties of the optimal solution and
proposed five heuristics to rapidly approximate the optimal solutions
of large-scale instances \citep{guan2022}. In all three papers \citep{Basciftci2020, Basciftci2021, guan2022},
detailed case studies with realistic data are conducted to demonstrate
the advantages and practicability of different methods under varied
circumstances; however, the possibilities of discovering the optimal
solution for large-scale ODMTS-DA problem still remains
elusive. Therefore, this paper aims to explore the availability of
modeling and computational methods that can both rapidly solve
normal-sized instances and obtain the optimal solution for large-scale
instances in a reasonable amount of time. It should be noted that the model presented in the work by \citet{Basciftci2021} (the second ODMTS-DA study) is referenced in the subsequent sections of this paper as the \textit{bilevel model}, serving as the established framework for defining the ODMTS-DA problem.

Generally speaking, the newly proposed path-based methods in this
paper intends to reformulate the bilevel formulation of the ODMTS-DA
problem into a single-level MIP. The reformulation not only provides
an alternative view of the original problem but can also ease the
computational difficulties. Multiple previous studies have applied
reformulations on bilevel optimization in the field of transportation
science and demonstrated their advantages, especially on enhancing
computational efficiency on large-scale instances
\citep{meng2001equivalent, marcotte2009toll, calvete2014planning,
  alizadeh2019dynamic}.  These reformulations mainly utilize
Karush–Kuhn–Tucker (KKT) conditions of the lower level problem or
strong-duality theorem.  For instance, \citet{goerigk2017line}
presented two different reformulation approaches to solve the line
planning with route choice problem.  The first reformulation utilizes
dualization of the inner-level routing problem, and the second one is
based on additional shortest-path constraints. In another example,
\citet{kara2004designing} first formulated the hazardous materials
transportation problem as a bilevel model; they then utilized
slackness conditions to reformulate the problem as a single-level MIP.

In summary, the proposed methods introduced in this paper are novel
from the following perspectives to address the ODMTS-DA problem: (i)
the paper proposes a path-based reformulation to develop
single-level optimization models satisfying the optimality of the
lower level problem, (ii) the reformulation only requires a small
number of critical paths, (iii) the approach generalizes previous
studies by considering a more flexible bilevel formulation and enabling choice function for riders as a black-box, and (iv)
the paper proposes preprocessing techniques that significantly reduce the
formulation sizes, resulting in single-level MIP reformulations that
can be directly solved by commercial solvers.


\section{The ODMTS Design with Adoptions (ODMTS-DA) Problem}
\label{sect:ProblemStatement}

This section introduces the main problem of this study---ODMTS Design
with Adoptions (ODMTS-DA). The goal for proposing this problem arises
from a desire to construct an ODMTS design $\mathbf{z}$ that
concurrently serves the existing transit users and captures latent
demand with the choices of riders. 
In general, the ODMTS-DA determines the locations (modeled by using arc variables $\mathbf{z}$) for the deployment of new fixed-route buses and provides multimodal pathways comprising fixed routes and on-demand shuttles to riders. The key assumption in ODMTS-DA is: it is framed around two distinct rider groups: (i) core riders, committed to utilizing the ODMTS, and (ii) latent riders, who have the option to choose between the ODMTS and other transportation modes. These mode choices are modeled using a predefined choice model that can be considered as a black-box. The primary aim of the ODMTS-DA problem is to minimize an aggregate weighted cost incurred by the system, stemming from three components: (i) operation of the newly-designed fixed routes, (ii) serving the core ridership, and (iii) serving the adopted latent demand.

The comprehensive details and rationale for the ODMTS-DA are elaborated in the subsequent subsections. Table~\ref{table:nomenclature}
outlines all nomenclature used in this paper, and they are
specifically discussed in Sections~\ref{sect:ProblemStatement},
\ref{sect:paths_method},  and~\ref{sect:computational_enhancements}. In this
section, Section~\ref{subsect:preliminaries} first introduces the
problem setting, and Section~\ref{subsect:bilevel} summarizes a
bilevel optimization framework that formally defines the ODMTS-DA
problem.

\begin{table}[!t]%
\SingleSpacedXI
    \centering%
    \resizebox{\textwidth}{!}{
    \begin{tabular}{l l }%
   
    \textbf{Symbol} & \textbf{Definition} \\
    \midrule
    \textbf{Sets}: &  \\
    $N$ & The set of ODMTS stops \\
    $H$ & The set of ODMTS hubs, $H \subseteq N$ \\
    $T$ & The set of all trips where each trip $r$ is defined by an unique O-D pair: ${or}^r, {de}^r \in N$ \\ 
    $T_{latent}$ & The set of latent trips, i.e., the latent demand \\
    $T_{core}$ & The set of core trips, i.e., existing transit demand \\
    $Z$& The set of hub arcs that are considered in the ODMTS-DA problem \\
    $Z_{fixed}$ & The set of hub arcs that are fixed in the ODMTS design $Z_{fixed} \subseteq Z$ \\
    $V^r, E^r$ & The sets of nodes or edges for a particular trip $r$ \\
    $L^r$ & The set of allowed on-demand shuttle arcs for trip $r$ \\
    $\Pi^r$ & The set of all valid ODMTS paths that connect ${or}^r$ and ${de}^r$ ($r \in T_{latent}$) \\
    $A^r$ & The set of ODMTS paths that a latent trip $r$ would \textbf{adopt} ($A^r \subseteq \Pi^r \enspace \forall r \in T_{latent}$) \\
    ${AP}^r$ & The set of \textbf{profitable} ODMTS paths that a latent trip $r$ would \textbf{adopt} \\
    ${ANP}^r$ & The set of \textbf{non-profitable} ODMTS paths that a latent trip $r$ would \textbf{adopt} \\
    ${RP}^r$ & The set of \textbf{profitable} ODMTS paths that a latent trip $r$ would \textbf{reject} (${RP}^r \subseteq \Pi^r \enspace \forall r \in T_{latent}$)  \\
    ${RNP}^r$ & The set of \textbf{non-profitable} ODMTS paths that a latent trip $r$ would \textbf{reject} \\ 
    $T_{latent,red}, L_{red}^r, Z_{red}^r$ & The \textbf{reduced} versions of the previously introduced  trip and arc sets after preprocessing \\ 
    $\Pi_{red}^r, A_{red}^r, {RP}_{red}^r$ & The \textbf{reduced} versions of the previously introduced  path sets after preprocessing \\
    $x(\pi)$ & The set of hub legs used by path $\pi$ for trip $r$ \\
    $y(\pi)$ & The set of shuttle legs used by path $\pi$ for trip $r$ \\
    ${\cal O}^r_f({\bf z})$ & The set of optimal paths for trip $r$ given a design $\mathbf{z}$ \\
    \midrule
    
    \textbf{Parameters}: & \\
    $\theta$ & Convex combination factor for weighted cost, $\theta \in [0, 1]$ \\
    $d_{ij}$ & Car travel distance between stops $i, j \in N$ \\
    $t_{ij}$ & Car travel time between stops $i, j \in N$ \\
    $d_{hl}'$ & Bus travel distance between hubs $h, l \in H$, not applicable to rail \\
    $t_{hl}'$ & Bus or rail travel time between hubs $h, l \in H$ \\
    $p^r$ & Number of passengers that take trip $r$ \\
    $n_{hl}$ & Number of buses or trains operating between hubs $h, l \in H$ during the planning horizon\\
    $t_{hl}^{wait}$ & A rider's expected waiting time for a bus or a train that operates between hubs $h,l \in H$ \\
    $b_{dist}, b_{time}$ & Operating cost for newly-designed bus arcs, with different measurements\\
    $\beta_{hl}$ & Weighted cost of investing a bus or train leg between hubs $h, l \in H$ \\
    $\tau_{hl}^r$ & A rider's weighted cost when traveling between hubs $h, l \in H$ with buses or rail for trip $r$\\
    $\gamma_{ij}^r$ & A rider's weighted cost when traveling between stops $i, j \in N$ with shuttles for trip $r$ \\
    $\omega_{dist}, \omega_{time}$ & Shuttle operating cost, with different measurements \\
    $\phi$ & Fixed ticket price charged by the transit agencies \\
    $\varphi$ & Transit agency's revenue for each rider in terms of weighted cost \\
    $\bar{g}^r$ & Upper bound for the weighted cost for a trip $r$\\
    $\underline{g}^r$ & Lower bound for the weighted cost for a trip $r$ \\ 
    $t_{\pi}^r$ & A rider's travel time for trip $r$ under the ODMTS path $\pi$\\
    $l_{\pi}^r$ & Number of transfers for trip $r$ under the ODMTS 
    path $\pi$\\
    $l_{ub}^r$ & Transfer tolerance for a latent trip $r \in T_{latent}$ in deciding ODMTS adoption, used in choice model~\eqref{eq:ChoiceFunctionTimeTransfer}\\
    $\alpha^r$ & Adoption factor of a latent trip $r \in T_{latent}$, used in choice models~\eqref{eq:ChoiceFunctionTime} and~\eqref{eq:ChoiceFunctionTimeTransfer}\\
    $t_{cur}^r$ & Travel time of a latent trip $r \in T_{latent}$ with its current travel mode, used in choice models~\eqref{eq:ChoiceFunctionTime} and~\eqref{eq:ChoiceFunctionTimeTransfer}\\

    \midrule
    \textbf{Decision Variables}: & \\
    $z_{hl}$   & A binary variable indicating if a hub arc between hubs $h, l \in H$ is open \\
    $x_{hl}^r$ & A binary variable indicating if a hub arc $(h, l)$ is used by trip $r$ \\
    $y_{ij}^r$ & A binary variable indicating if a shuttle arc $(i, j)$ is used by trip $r$ \\
    $t^r$ & Travel time of trip $r$ for a rider, dynamically based on ODMTS paths  \\
    $g^r$ & Weighted cost of trip $r$ for a rider, dynamically based on ODMTS paths \\
    $m^r$ & Minimal weighted cost of trip $r$ for a rider, dynamically based on ODMTS design $\mathbf{z}$ \\
    $\delta^r$ & A binary variable indicating if latent trip $r \in T_{latent}$ adopts \\
    $f_\pi^r$ & A binary variable indicating if path $\pi$ is \textbf{feasible} for latent trip $r$, dynamically based on  $\mathbf{z}$  \\
    $\lambda_{\pi}^r$ & A binary variable indicating if path $\pi$ is assigned to latent trip $r$ \\
    
    \bottomrule
    \end{tabular}
    }
    \caption{Nomenclature used in this study. For simplicity, $\pi^r$ is simplified as $\pi$ in variables and parameters.}
    \label{table:nomenclature}
\end{table}

\subsection{ODMTS Preliminaries}
\label{subsect:preliminaries}

The ODMTS considers a set of nodes $N$, that represents the stops for
on-demand shuttles to pick-up or drop-off riders.  Given a pair of
stops $i, j \in N$, the direct travel time and travel distance between
these stops are represented by $t_{ij}$ and $d_{ij}$, respectively,
which correspond to the travel time and distance for on-demand
shuttles serving these stops. A subset $H \subseteq N$, which
corresponds to the hubs, can be employed to establish high-frequency
buses or rails. The decision on hub locations must be made prior to any computational processes, usually involving the selection of a few key locations to ensure comprehensive coverage throughout the service area. There are two types of hub-to-hub connections
considered in the ODMTS. The first type is defined as
\textit{newly-designed bus arcs} or simply \textit{new bus arcs}, where their connections are decided by the optimization models. Note that
the ODMTS normally does not design any new rail arcs due to practical
considerations. The second type is referred as \textit{backbone arcs},
they represent pre-selected existing transit lines such as rails and
Bus Rapid Transit (BRT), and can be integrated into the ODMTS as fixed
arcs, i.e., not decided by the optimization models. Similarly, for
buses and rails operating between hubs $h, l \in H$, $t_{hl}'$ and
$d_{hl}'$ are defined to represent time and distance values. The inclusion of $t_{hl}'$ and $d_{hl}'$ for hub-to-hub connections is necessitated by the potential variation in travel time and distance between the same origin-destination pair for different modes.

The set $T$ includes all trips that are considered by the ODMTS,
including both the existing and latent demand. Three values are
associated with each trip $r \in T$: (i) an origin stop $or^r \in N$,
(ii) a destination stop $de^r \in N$, and (iii) the number of riders
taking that trip $p^r \in \mathbb{Z}_{+}$. In order to serve a
particular trip $r$, an ODMTS path $\pi$ that connects ${or}^r$ and
${de}^r$ is proposed by the transit agency, and it is usually
multimodal. For example, a typical ODMTS path consists of three
trip-legs: a shuttle leg that serves the first-mile, a hub leg
in the middle (taking rail or bus), and another shuttle leg that serves the last-mile. A
path is allowed to include multiple hub legs, while
shuttle legs can only be employed to serve as the first or last
legs. In this study, the term \textit{direct shuttle path}, indicates an ODMTS path that is covered by a single-leg
on-demand shuttle that connects ${or}^r$ and ${de}^r$. 
In summary, on-demand shuttles are only allowed for three types of trips: (i) the direct trip, (ii) a
first leg from trip origin to a hub, and (iii) a last leg from a hub to trip destination. In addition, note that paths solely served with two on-demand shuttle legs are not considered. For example, a path in the form ${origin} \rightarrow {hub} \rightarrow {destination}$ with two shuttles are not allowed.


To have a weighted objective considering both convenience and cost
aspects of the ODMTS design, a parameter $\theta \in [0,1]$ is
employed such that the trip duration associated with convenience are
multiplied by $\theta$ and expenses are multiplied by $1 - \theta$. When a new bus arc connects hubs $h, l \in H$, $n_{hl}$ denotes the number of buses operating along this route within the planning horizon, where $n_{hl}$ is an input parameter defined by the transit agency for setting the hourly frequency. For instance, setting $n_{hl} = 12$ in a four-hour time horizon implies a bus departing from hub $h$ to hub $l$ every 20 minutes $(4 \times 60 / 12)$. Alternatively, the agency could also specify an hourly frequency of 3 for the bus service between $h$ and $l$, $n_{hl} = 3 \times 4$. It follows that a weighted cost $\beta_{hl}$ is
considered by the transit agency using equation $\beta_{hl} = (1 -
\theta) n_{hl} \, d_{hl}' \, b_{dist}$, where $b_{dist}$ is the cost
of operating a bus per kilometer and $d_{hl}'$ is the bus travel
distance from $h$ to $l$. It is worth pointing out that $\beta_{hl}$
can be alternatively modeled by $\beta_{hl} = (1 - \theta) n_{hl} \,
t_{hl}' \, b_{time}$, where $b_{time}$ and $t_{hl}'$ represent the
cost of operating a bus per hour and bus travel time from hub $h$ to
$l$, respectively. Furthermore, for backbone hub connections between
$h$ and $l$, e.g., a rail arc that is a part of the existing rail
system, $\beta_{h,l} = 0$ is applied instead of the abovementioned
equations.

In this paper, maintaining $n_{hl}$ as a constant for all hub arcs is required. However, as noted by \citet{dalmeijer2020transfer}, to give more flexibility to the network designers, the ODMTS can be extended to accommodate multiple predefined hourly frequency choices for each hub arc. Consequently, the optimization model could select the most suitable option from these choices. Since this study considers a complex transit network design problem by integrating rider choice decisions within a bilevel optimization framework, the decisions involving the selection of hourly frequency for hub arcs, shuttle ride-sharing, shuttle flow reservations, and shuttle fleet sizing decisions of on-demand shuttles are omitted while designing the network. These decisions can be either addressed by modifying constraints or considering them in the subsequent levels of decision-making, which involve finer level operational planning under a given network design.


In addition to the weighted cost of opening legs between the hubs,
serving each trip has its weighted cost to the transit agency. This
cost consists of two parts depending on the path utilized by the trip,
which are hub legs and shuttle legs. For the former one, as
operating cost of the open legs between hubs are considered within the
investment, only inconvenience cost of each trip $r \in T$, denoted as
$\tau_{hl}^r$, is computed by $\tau_{hl}^r = \theta (t_{hl}' +
t_{hl}^{wait}) $, where $t_{hl}^{wait}$ is the expected waiting time
of a rider for the bus or train between hubs $h$ and $l$. On the other
hand, for the parts of the trips traveled by on-demand shuttles, the
transit agency incurs $\gamma_{ij}^r$ for each trip $r \in T$ using
equation $\gamma_{ij}^r = (1 - \theta) \omega_{dist} \, d_{ij} + \theta
t_{ij} $, where $\omega_{dist}$ is the cost of operating a shuttle per
kilometer. Moreover, $\gamma_{ij}^r$ can be modeled alternatively by using $\gamma_{ij}^r = (1 - \theta) \omega_{time} \, t_{ij} + \theta
t_{ij} $, where $\omega_{time}$ is the shuttle operating cost per hour.


\subsection{Problem Description and Bilevel Formulation}
\label{subsect:bilevel}

This section formally defines the ODMTS-DA problem as a bilevel
optimization problem. The formulation generalizes the
proposal by \citet{Basciftci2021} to allow for more flexibility in the
mode choice models. Figure~\ref{fig:bilevel} presents the bilevel
optimization framework for ODMTS-DA problem.

\paragraph{Decision Variables}
Binary variable $z_{hl}$ represents whether a hub arc between $h,l \in H$ is open. Backbone arcs form a set $Z_{fixed}
\subseteq \{ (h,l) : \forall h, l \in H \}$, and their corresponding
$z_{hl}$ are always open. For each trip $r \in T$, binary variables
$x_{hl}^r$ and $y_{ij}^r$ indicate whether trip $r$ takes the hub leg
between $h,l \in H$, and the shuttle leg between stops $i, j
\in N$, respectively. An ODMTS path $\pi$ in the network for trip $r$ is defined as a sequence of distinct opened hub arcs and shuttle arcs that connect a sequence of vertices starting from origin ${or}^r$ and ending at destination ${de}^r$. Under a given solution, a path of trip $r$ can be constructed from the open arcs specified by the $x_{hl}^r$
and $y_{ij}^r$ variables.

\paragraph{Adoption Choices}
To integrate the adoption choices of riders into the formulation, the
trip set $T$ is divided into two subsets with respect to the riders'
adoption characteristics: (i) a \textit{core trip set} $T \setminus
T_{latent}$ representing the existing demand and (ii) a \textit{latent trip
  set} $T_{latent}$ representing the latent demand. Core trips correspond to
the set of riders of the existing transit system: they are assumed to
continue adopting the ODMTS.  Latent trips correspond to the set of riders who are currently
traveling with other modes and might switch to transit due to the
deployment of ODMTS. The transit agency charges each rider a ticket
fare $\phi$ to use the ODMTS, irrespective of the assigned
paths. Hence, a fixed value $\varphi$, computed by $\varphi = (1 -
\theta)\phi$, becomes an additional weighted revenue to the transit
agency for riders switching to the ODMTS.

In this study, the adoption behavior of a latent trip $r \in T_{latent}$ is
assumed to depend on the features of an assigned path $\pi$, such as
trip duration $t_\pi^r$ and the number of transfers $l_\pi^r$, and
these features can be computed by the the following equations:

\begin{equation} 
t_\pi^r = \sum_{h,l \in H}  (t_{hl}' + t_{hl}^{wait}) x_{hl}^r + \sum_{i,j \in N}  t_{ij} y_{ij}^r
\end{equation}
\begin{equation} 
l_\pi^r = \sum_{h,l \in H} x_{hl}^r + \sum_{i,j \in N} y_{ij}^r - 1
\end{equation}

\noindent
For each latent trip $r$, a choice model ${\cal C}^r$ returns a binary
adoption decision $\delta^r$, which is equal to 1 when its riders
decide to adopt the ODMTS given the proposed path $\pi$.  The
structure of ${\cal C}^r$ can vary from a threshold model to a more
complicated machine learning model. The choice model can be
  arbitrarily complicated in this paper since it is abstracted in the
  path enumeration. \cite{Basciftci2021} proposed the choice model
\begin{equation} 
\label{eq:ChoiceFunctionTime}
{\cal C}^r(\mathbf{x}^r, \mathbf{y}^r) = {\cal C}^r(\pi) \equiv \mathbbm{1}(t_\pi^r \leq \alpha^r \ t_{cur}^r) \quad \forall r \in T_{latent}
\end{equation}
that directly compares the ODMTS trip duration $t_\pi^r$ with the
travel duration $t^r_{cur}$ under the current mode multiplied by a
parameter $\alpha^r$. In the computational experiments, this paper
also explores the following choice model
\begin{equation} 
    \label{eq:ChoiceFunctionTimeTransfer}
    {\cal C}^r(\mathbf{x}^r, \mathbf{y}^r) = {\cal C}^r(\pi) \equiv \mathbbm{1}(t^r_{\pi} \leq \alpha^r \ t^r_{cur}) \land \mathbbm{1}(l^r_{\pi} \leq l^r_{ub}) \quad \forall r \in T_{latent}
\end{equation}
which takes transfers into account, i.e., riders reject the ODMTS
when they are assigned to a path with too many transfers. For each $r
\in T_{latent}$, the transfer tolerance is denoted as $l_{ub}^r$ and $l_\pi^r$
represents the number of transfers in a path $\pi$. For example, consider a path $\pi$ with 4 transfers $(l_\pi^r = 4)$ offered to riders in trip $r$: $shuttle \rightarrow bus \rightarrow rail \rightarrow bus \rightarrow shuttle$. If the riders in $r$ have a transfer tolerance of 2 $(l_{ub}^r = 2)$, then according to choice model \eqref{eq:ChoiceFunctionTimeTransfer}, the riders in $r$ will reject the proposed path $\pi$.

\paragraph{The Bilevel Optimization Model}

\begin{figure}[!t]
  \begin{subequations}
\label{eq:bl}    
\SingleSpacedXI
\begin{alignat}{1}
\min_{z_{hl}} \quad & \sum_{h,l \in H} \beta_{hl} z_{hl} + \sum_{r \in T_{core}} p^r g^r + \sum_{r \in T_{latent}} p^r \delta^r (g^r - \varphi) \label{eq_bi:upperLevelObj} \\
\text{s.t.} \quad 
& \sum_{l \in H} z_{hl} = \sum_{l \in H} z_{lh} \quad \forall h \in H \label{eq_bi:upperLevelConstr1} \\
& z_{hl} = 1 \quad \forall (h, l) \in Z_{fixed} \label{eq_bi:fixed_arc} \\
& \delta^r = {\cal C}^r(\mathbf{x^r}, \mathbf{y^r}) \quad \forall r \in T_{latent} \label{eq_bi:userChoiceModel} \\
& z_{hl} \in \{0,1\} \quad \forall h,l \in H  \label{eq_bi:binaryConstraint} \\
& \delta^r \in \{0,1\} \quad \forall r \in T_{latent} \label{eq_bi:continuousConstraint} 
\end{alignat}
\end{subequations}

where $(\mathbf{x^r}, \mathbf{y^r}, g^r)$ are a solution to the optimization problem

\begin{subequations}
\label{eq_bi:lowerLevelProblem}
\SingleSpacedXI
\begin{alignat}{1}
(\mathbf{x^r}, \mathbf{y^r}, g^r) \in \argmin_{x_{hl}^r, y_{ij}^r} \quad &  g^r = \sum_{h,l \in H} \tau_{hl}^r x_{hl}^r + \sum_{i,j \in N} \gamma_{ij}^r y_{ij}^r \label{eq_bi:lowerLevelObj} \\
\text{s.t.} \quad
& \sum_{\substack{h \in H \\ \text{if } i \in H}} (x_{ih}^r - x_{hi}^r) + \sum_{i,j \in N}  (y_{ij}^r - y_{ji}^r) =             \begin{cases}
     1 & \text{if  } i = or^r \\
    -1 & \text{if  } i = de^r \\
    0 & \text{otherwise}
    \end{cases} 
\quad \forall i \in N \label{eq_bi:minFlowConstraint} \\
& x_{hl}^r \leq z_{hl} \quad \forall h,l \in H \label{eq_bi:openFacilityOnlyAvailable} \\
& x_{hl}^r, y_{ij}^r \in \{0,1\} \quad \forall i,j \in N \quad \forall h,l \in H \quad \forall r \in T \label{eq_bi:integralityFlowConstr}
\end{alignat}
\end{subequations}
\caption{The Bilevel Optimization Model for the ODMTS Design with Adoption (ODMTS-DA) Problem}
\label{fig:bilevel}
\end{figure}

\begin{figure}[!t]
\begin{subequations}
\label{eq_bi:lowerLevelProblem_sc}
\SingleSpacedXI
\begin{alignat}{1}
\lexmin_{x_{hl}^r, y_{ij}^r, g^r, t^r} \quad &  \langle g^r, t^r \rangle \label{eq_bi:lowerLevelObj_sc} \\
\text{s.t.} \quad
& g^r = \sum_{h,l \in H} \tau_{hl}^r x_{hl}^r + \sum_{i,j \in N} \gamma_{ij}^r y_{ij}^r \label{eq_bi:drDefinition_sc} \\
& t^r = \sum_{h,l \in H}  (t_{hl}' + t_{hl}^{wait}) x_{hl}^r + \sum_{i,j \in N}  t_{ij} y_{ij}^r \label{eq_bi:frDefinition_sc} \\
& \eqref{eq_bi:minFlowConstraint}-\eqref{eq_bi:integralityFlowConstr} \notag
\end{alignat}
\end{subequations}
\caption{The Subproblem with a Unique Optimal $(g^r, t^r)$ Solution Proposed in \citet{Basciftci2021}.}
\label{fig:original_subproblem}
\end{figure}

The leader problem (Equations
\eqref{eq_bi:upperLevelObj}--\eqref{eq_bi:continuousConstraint})
designs the network between the hubs for the ODMTS, whereas the
follower problem (Equations
\eqref{eq_bi:lowerLevelObj}--\eqref{eq_bi:integralityFlowConstr})
obtains paths for each trip $r \in T$ by considering the legs of the
network design and the shuttles to serve the first and last miles.
The objective of the leader problem minimizes the sum of (i) the
weighted cost of opening new bus legs, (ii) the weighted cost of the
core trips, and (iii) the weighted net cost of the riders with choice
who are switching to the ODMTS.  
Note that the third component of the objective function includes a nonlinear term as the binary variable $\delta^r$ is utilized to represent the riders that adopt ODMTS. Constraint~\eqref{eq_bi:upperLevelConstr1} ensures weak connectivity between the
hubs such that each hub has the same number of incoming and outgoing
open arcs. Constraint~\eqref{eq_bi:fixed_arc} enforces the fixed arcs to be opened. Constraint~\eqref{eq_bi:userChoiceModel} represents the
adoption decision of the riders in $T_{latent}$ depending on the ODMTS path
constructed by $\mathbf{x}^r$ and $\mathbf{y}^r$.

For each trip $r$, the objective of the follower problem
\eqref{eq_bi:lowerLevelProblem} minimizes the cost and inconvenience (weighted cost), with the solution of the leader problem ($\mathbf{z}$) serving as inputs. Constraint~\eqref{eq_bi:minFlowConstraint} guarantees
flow conservation for the bus, rail, and shuttle legs used in trip $r$
for origin, destination, and each of the intermediate
points. Constraint~\eqref{eq_bi:openFacilityOnlyAvailable} ensures
that a path only uses open legs between the hubs. Given that the
follower problem has a totally unimodular constraint matrix, it can be
solved as a linear program.

The bilevel optimization model can be understood intuitively as
follows.  The transit agency obtains an ODMTS design by minimizing the
weighted costs of the design and the trips considered. The paths derived from this
design are proposed to riders, and riders with choice determine
whether they adopt the ODMTS based on their suggested paths or stay
with their previous modes. It is worth noticing that the paths are proposed by the agency instead of the riders, this is the key difference between ODMTS and regular transit systems. The subproblem guarantees that the proposed
paths are minimizing the cost and inconvenience of the riders,
preventing the formulation to propose suboptimal paths that riders
would systematically reject. The optimal design can thus be viewed as
an equilibrium point for such a game between the transit agency and
the latent demand.

\paragraph{The Original Bilevel Optimization Model}

The follower subproblem in~\eqref{eq_bi:lowerLevelProblem} may have
multiple optimal solutions that may lead to routes with different
travel durations. \citet{Basciftci2021} eliminates this issue by
proposing the subproblem depicted in Figure
\ref{fig:original_subproblem}.  For each trip $r$, the objective of
the follower problem~\eqref{eq_bi:lowerLevelProblem_sc} minimizes the
lexicographic objective function $\langle g^r, t^r \rangle$, where
$g^r$ corresponds to the cost and inconvenience of trip $r$ and $t^r$
breaks potential ties between the solutions with the same value of
$g^r$ by returning the most convenient path for the rider of trip
$r$. \citet{Basciftci2021} broke ties this way because it aligned with
the choice model they used (i.e.~\eqref{eq:ChoiceFunctionTime}), that is solely based on travel
duration. They also showed that the lexicographic minimizer of problem
\eqref{eq_bi:lowerLevelProblem_sc} exists and results in a unique $(g^r, t^r)$ solution. This unicity
property simplifies the design of algorithms and is exploited by
\citet{Basciftci2021}.

To solve this resulting bilevel optimization problem,
\citet{Basciftci2021} proposed an exact decomposition algorithm that
combines Benders decomposition algorithm with a combinatorial cut
generation procedure to integrate rider adoption constraints. Despite
the addition of valid inequalities and the application of several
preprocessing techniques, the computational studies revealed some of
the limitations of this approach for large-scale instances (e.g.,
metropolitan cities): those instances have large optimality gaps and
run times.

The formulations and algorithms proposed in this study do not require
this unicity property for the subproblem. This provides a more
flexible framework where the subproblems do not need to have unique
optimal solutions. The experimental results will consider both the
general formulation and the more specialized model when comparing
execution times.

\section{Path-Based Formulations}
\label{sect:paths_method}

This section presents path-based formulations for addressing the
ODMTS-DA problem. Section~\ref{sect:paths} discusses the nature of the paths and 
the set of paths useful for the path-based formulations along with their properties.
Sections~\ref{subsect:C-Path} and~\ref{subsect:P-Path} present the two path-based formulations. 

\subsection{Paths and Their Properties}
\label{sect:paths}

\paragraph{Paths} Given a design ${\bf z}$, the follower subproblem
\eqref{eq_bi:lowerLevelProblem} returns a set of paths for each rider.
For trip $r$, these paths are specified by the vectors ${\bf x^r}$
and ${\bf y^r}$, i.e., they specify which hub legs and which shuttle
legs comprise the path. In fact, given the restrictions on shuttle legs (which
can only be the first and last legs of a path or a direct path), a path is
uniquely specified by its set of hub legs. This property is used by
the path formulations.

The path formulations reason directly in terms of paths, in addition
to the decision variables ${\bf x^r}$ and ${\bf y^r}$. For a rider
$r$, a path $\pi$ specifies which hub legs and shuttle legs trip $r$
uses if the path is adopted.  The notations
\begin{subequations}
  \begin{alignat*}{1}
 x(\pi) = \{ (h,l) \in \pi : h, l \in H \}  \\
 y(\pi) = \{ (i,j) \in \pi : i, j \in N \} 
\end{alignat*}
\end{subequations}
denote the set of hub legs and shuttle trips used by path $\pi$
respectively. The paper also uses ${\cal C}^r(\pi)$ to denote whether riders in
trip $r$ adopt path $\pi$.

\paragraph{Set of Paths} The path formulations use some specific sets of paths to
solve the
ODMTS-DA problem. These sets are only constructed for riders
in $T_{latent}$.  The set of paths for trip $r \in T_{latent}$ are constructed from a
graph $G^r = (V^r, E^r)$ that contains all the hub legs and the allowed shuttle
trips, where the set of nodes is defined as $V^r = \{or^r, de^r \} \cup H$,  the set of arcs as $E^r = Z \cup L^r$, with $L^r$ is the set of allowed shuttle trips, i.e., 
\begin{equation*} 
\label{eq:gen_set_allow_shuttle_arcs}
L^r = 
    \{  ({or}^r, {de}^r) \} 
    \cup 
    \{ ({or}^r, h):\forall h \in H \} 
    \cup 
    \{ (h, {de}^r):\forall h \in H \}  
\end{equation*}

The path-based formulations rely on three sets:
\begin{subequations}
  \begin{alignat*}{1}
 &   \Pi^r = \mbox{ the set of all paths in }G^r\mbox{ connecting }or^r\mbox{ to }de^r. \\
 &   A^r = \{ \pi \in \Pi^r : {\cal C}^r(\pi) = 1 \}. \\
 &   RP^r = \{ \pi \in \Pi^r : {\cal C}^r(\pi) = 0 \; \land \; g(\pi) < \varphi \}. 
\end{alignat*}
\end{subequations}
where 
\[
g(\pi) = \sum_{(h,l) \in x(\pi)} \tau_{hl} + \sum_{(i,j) \in y(\pi)}  \gamma_{ij}.
\]
The set $A^r$ represents the set of all paths that trip $r$ would
adopt. The set ${RP}^r$ denotes the set of all paths that trip $r$
would reject but which profitable from the perspective of the transit
agency in terms of weighted cost. 
Algorithm~\mynameref{alg:s1_alg_general_form} in Appendix~\ref{sect:pathEnumeration} can be used to 
obtain these sets. In practice, these path sets are generated using various pre-processing techniques for computational efficiency, as discussed in Section~\ref{sect:computational_enhancements}.  Why these sets are useful will become clear
subsequently. The first path-based formulation uses sets $\Pi^r$ and
$A^r$ for each latent trip $r \in T_{latent}$, whereas the second formulation is
based on ${RP}^r$ and $A^r$. In general, ${RP}^r$ is
smaller than $\Pi^r$ because of the rare appearances of profitable paths.  For proving the results, it is also useful to
define the set
\begin{subequations}
  \begin{alignat*}{1}
 &   {RNP}^r = \{ \pi \in \Pi^r : {\cal C}^r(\pi) = 0 \; \land \; g(\pi) \geq \varphi \}. 
\end{alignat*}
\end{subequations}
Note that $A^r$, ${RP}^r$, and ${RNP}^r$ form a partition of $\Pi^r$. It is also
useful to further partition $A^r$ into ${AP}^r$ and ${ANP}^r$ with 
\begin{subequations}
  \begin{alignat*}{1}
    &   {AP}^r = \{ \pi \in A^r : g(\pi) < \varphi \}, \\
    &   {ANP}^r = \{ \pi \in A^r : g(\pi) \geq \varphi \}.   
\end{alignat*}
\end{subequations}

\paragraph{Path Formulations}
It is important to recall that the bilevel nature of the ODMTS-DA
problem is only due to the riders with choice: the bilevel formulation
could include equations~\eqref{eq_bi:minFlowConstraint} and
\eqref{eq_bi:openFacilityOnlyAvailable} for riders in $T_{core}$
in the leader problem. Indeed, since riders in $T_{core}$ must
use the ODMTS, the optimization model will necessarily minimize their
objective terms. This is not the case of riders of choice, since their
choice function decides whether they adopt the ODMTS. 
Without a
bilevel model, the single-level optimization, which directly incorporates the lower level constraints to the upper level problem for riders in $T_{latent}$, can intentionally propose paths to riders with choices that they will reject, if no profitable paths
exist for them. This happens as the optimality of the lower level is not ensured in this single-level problem for these riders, and the suggested routes and choices of the riders are evaluated only from the perspective of the transit agency. 
This is discussed at length by
\cite{Basciftci2021}, where the authors studied this single-level variant of the problem for comparison. Their computational study presents that this
single-level formulation only evaluates the suggested routes and choices of the riders from the
perspective of the transit agency, who consequently can suggest longer routes to the riders with
choice if serving them is not profitable. Thus, their inconvenience is explicitly omitted in the
system, which is undesirable for ensuring the access to the transit system. Consequently, the
authors observe significantly less riders adopting ODMTS under this formulation, in comparison
to the bilevel model. On the other hand, the bilevel formulation eliminates this pathological and
unfair behavior. This is aligned with the objectives of many transit agencies which aims at using
ODMTS to improve mobility for underserved communities.

As a result, in any single-level reformulation and solution algorithm, the bilevel nature of the
problem and the combinatorial choice functions of the riders need to be carefully addressed.
To this end, in this paper, both path formulations only reason
about paths for riders with choices: they continue to use decision
variables ${\bf x^r}$ and ${\bf y^r}$ for riders in $T_{core}$.

\paragraph{Path Properties} It is also useful to characterize the behavior
of the bilevel problem~\eqref{eq:bl} and the subproblem
\eqref{eq_bi:lowerLevelProblem}. Given a design ${\bf \bar{z}}$ and a
trip $r \in T_{latent}$, the subproblem returns a set of optimal paths ${\cal
  O}^r_f({\bf \bar{z}})$. The following two properties of the bilevel
problem~\eqref{eq:bl} are important in case ${\cal O}^r_f({\bf
  \bar{z}})$ is not a singleton, i.e., multiple paths have the same $g(\pi)$ value.

\begin{properties}
\label{prop:selectA}
  Given a design ${\bf \bar{z}}$, if there exists $\pi \in AP^r \cap
  {\cal O}^r_f({\bf \bar{z}})$, then the optimization selects $\pi$
  since it decreases the objective~\eqref{eq_bi:upperLevelObj} (the third term in the objective function is negative).
\end{properties}

\begin{properties}
\label{prop:selectRNP}  
Given a design ${\bf \bar{z}}$, if $AP^r \cap
  {\cal O}^r_f({\bf \bar{z}}) = \emptyset$ and 
there exists $\pi \in RNP^r \cap {\cal
  O}^r_f({\bf \bar{z}})$, then the optimization selects $\pi$,
since it does not increase the objective value.
\end{properties}

\noindent
Property~\ref{prop:selectA} favors the selection of a path in $AP^r$
over a path in $RP^r$, while Property~\ref{prop:selectRNP} prefers a
path in $RNP^r$ over a path in $ANP^r$, since the latter would induce
a positive weighted cost in the third term of the objective function.

\subsection{Formulation C-Path}
\label{subsect:C-Path}

This section introduces the first path-based formulation, {\sc C-Path}, that
reasons over the sets $\Pi^r$ and $A^r$ for each trip $r \in T_{latent}$. By
reasoning about these paths, and not over the variables ${\bf x^r}$
and ${\bf y^r}$ for riders in $T_{core}$, {\sc C-Path} can be
expressed as a single-level formulation. The key to {\sc C-Path} is to
make sure that only optimal paths (i.e., those returned by the
follower subproblem in the bilevel formulation) are selected for
riders with choice.  A high-level presentation of {\sc C-Path} is
presented in Figure~\ref{fig:c-path-high-level}. 
The formulation uses two abbreviations
\[
{\sc Feasible}(\pi,{\bf z}) \equiv \forall (h,l) \in x(\pi): z_{hl} = 1
\]
and 
\[
{\sc Optimal}(\pi,{\bf z}) \equiv \mathbbm{1}(g(\pi) = \min \{ g(\pi'): \pi' \in \Pi^r \; \land \; {\sc Feasible}(\pi',{\bf z}) \})
\]
${\sc Feasible}(\pi,{\bf z})$ holds if path $\pi$ is feasible under
design ${\bf z}$. Note that there is no need to consider the shuttle
  arcs since they are always available. ${\sc Optimal}(\pi,{\bf z})$
holds if path $\pi$ has an optimal objective value among all the
feasible paths of design ${\bf z}$.

\begin{figure}[!t]
\begin{subequations}
\label{eq_cpath_hl}
\begin{alignat}{1}
\min_{z_{hl}} \quad & \sum_{h,l \in H} \beta_{hl} z_{hl} + \sum_{r \in T_{core}} p^r g^r + \sum_{r \in T_{latent}} \sum_{\pi \in A^r} p^r \lambda^r_\pi (g(\pi)  - \varphi) 
\label{eq_cpath_hl:obj_func} \\
\text{s.t.} \quad & \sum_{l \in H} z_{hl} = \sum_{l \in H} z_{lh} \quad \forall h \in H  
\label{eq_cpath_hl:bus_flow_balance} \\
& z_{hl} = 1 \quad \forall (h, l) \in Z_{fixed} \label{eq_cpath_hl:fixed_arc} \\
& g^r = \sum_{h,l \in H} \tau_{hl}^r x_{hl}^r + \sum_{i,j \in N} \gamma_{ij}^r y_{ij}^r \quad \forall r \in T_{core} \label{eq_cpath_hl:grDefinition} \\
& \sum_{\substack{h \in H \\ \text{if } i \in H}} (x_{ih}^r - x_{hi}^r) + \sum_{i,j \in N}  (y_{ij}^r - y_{ji}^r) = \begin{cases}
 1 & \text{if  } i = or^r \\
-1 & \text{if  } i = de^r \\
0 & \text{otherwise}
\end{cases} \quad \forall i \in N, r \in T_{core} \label{eq_cpath_hl:minFlowConstraint} \\
& x_{hl}^r \leq z_{hl} \quad \forall h,l \in H, r \in T_{core} \label{eq_cpath_hl:openFacilityOnlyAvailable} \\
& \lambda_{\pi}^r \leq {\sc Feasible}(\pi,{\bf z}) \wedge {\sc Optimal}(\pi,{\bf z}) \quad \forall r \in T_{latent}, \pi \in \Pi^r \label{eq_cpath_hl:optimality} \\
& \sum\limits_{\pi \in \Pi^r} \lambda_\pi^r = 1 \quad \forall r \in T_{latent} \label{eq_cpath_hl:must_take_1_path} \\
& z_{hl}, x_{hl}^r, y_{ij}^r \in \{0, 1\} \quad \forall h,l \in H, i,j \in N, r \in T_{core} \label{eq_cpath_hl:binaryConstraint_start} \\
&  \lambda_{\pi}^r \in \{0, 1\} \quad  \forall r \in T_{latent}, \pi^r \in \Pi^r \label{eq_cpath_hl:binaryConstraint_end} 
\end{alignat}
\end{subequations}
\caption{The High-Level Description of Formulation {\tt C-Path}.}
\label{fig:c-path-high-level}
\end{figure}

In formulation~\eqref{eq_cpath_hl}, variables $z_{hl}$ and variables
$x^r_{hl}$, $y^r_{ij}$ , and $g^r$ for $r \in T_{core}$ are
directly adopted from the bilevel formulation~\eqref{eq:bl}. The
constraints for riders with choice are expressed in terms of paths.
Variable $\lambda_\pi^r$ indicates whether path $\pi$ is selected
among the optimal solutions of the follower problem
\eqref{eq_bi:lowerLevelProblem} for trip $r \in T_{latent}$. Constraints~\eqref{eq_cpath_hl:must_take_1_path} ensure that only one such path
is selected. Constraints~\eqref{eq_cpath_hl:optimality} guarantee that
the selected paths are indeed feasible and optimal for the follower
subproblem. The first two terms of the objective function are the same
as in the bilevel formulation. The third term is more interesting.
The selected path $\pi$ for a trip $r \in T_{latent}$ only contributes to the
objective if it belongs to $A^r$, i.e., it is adopted by trip $r$.
This achieves the same effect as the choice function in the bilevel
formulation. It is interesting to note that formulation
\eqref{eq_cpath_hl} is single-level: the use of paths has eliminated
the need for follower problems.

It remains to linearize the two abbreviations. The linearization of
${\sc Feasible}(\pi,{\bf z})$ introduces a variable $f^r_\pi$ for each
trip $r \in T_{latent}$ and each path $\pi \in \Pi^r$ and the following constraints:
\begin{subequations}
  \begin{alignat*}{1}
 &   f_\pi^r \leq z_{hl} \quad \forall r \in T_{latent}, \pi \in \Pi^r, (h,l) \in x(\pi) \\
 &   {\sc valid}(\pi,{\bf z}) \leq f^r_{\pi} \quad \forall r \in T_{latent}, \pi \in \Pi^r
\end{alignat*}
\end{subequations}
where
\[
  {\sc valid}(\pi,{\bf z}) \equiv \sum\limits_{(h,l) \in x(\pi)} z_{hl} - |x(\pi)| + 1
  \]
This last constraint expresses that the path $\pi$ is allowed by
design ${\bf z}$. In other words, a path $\pi$ cannot be chosen if it is infeasible ($f_\pi^r = 0$) under design $
\mathbf{z}$.  The linearization of ${\sc Optimal}(\pi,{\bf z})$
introduces a variable $m^r$ for all $r \in T_{latent}$ and the constraints
\[
f^r_\pi \rightarrow m^r \leq g(\pi) \quad \forall r \in T_{latent}, \pi \in \Pi^r
\]
where $m^r$ represents the optimal objective value of the subproblem. To ensure that the selected
paths are feasible, the linearization of~\eqref{eq_cpath_hl:optimality} adds the constraints
\[
\lambda_\pi^r \leq f_\pi^r \quad \forall r \in T_{latent}, \pi \in \Pi^r.
\]
To guarantee that the selected paths are optimal, the linearization adds the constraints
\[
\lambda^r_\pi \rightarrow m^r \geq g(\pi) \quad \forall r \in T_{latent}, \pi \in \Pi^r.
\]
Indeed, if $m^r < g(\pi)$, then it must be the case that
$\lambda^r_\pi = 0$ and path $\pi$ cannot be selected. 

\begin{thm} \label{thm:CPathOptimality}
Formulation~\eqref{eq_cpath_hl} returns a design that is optimal for Formulation~\eqref{eq:bl}.
\end{thm}

Detailed proof of {\sc C-PATH}'s optimality and its full formulation can be found in Appendix~\ref{sect:c_path_additional}. The presentation of {\sc C-PATH}, the abbreviations used in  {\sc C-PATH}, and the constraints lineralization together serve as a transitional step, linking the bilevel formulation to the {\sc P-PATH} formulation, the main contribution of this paper. The similarity between {\sc C-PATH} and the bilevel formulation is that $({\bf x^r},{\bf y^r})$ are used for core trips in $T_{core}$, while {\sc P-PATH} deviates from this approach.

\subsection{Formulation P-Path}
\label{subsect:P-Path}

\begin{figure}[!t]
\begin{subequations}
\label{eq:ppathhl}
\begin{alignat}{1}
\min_{z_{hl}} \quad & \sum_{hl \in H} \beta_{hl} z_{hl} + \sum_{r \in T_{core}} p^r g^r + \sum_{r \in T_{latent}} \sum_{\pi \in A^r} p^r \lambda^r_\pi (g(\pi)  - \varphi) \\
\text{s.t.} \quad 
& \sum_{l \in H} z_{hl} = \sum_{l \in H} z_{lh} \quad \forall h \in H  \label{eq_ppath_hl:bus_flow_balance} \\
& z_{hl} = 1 \quad \forall (h, l) \in Z_{fixed} \label{eq_ppath_hl:fixed_arc} \\
& g^r = \sum_{hl \in H} \tau_{hl}^r x_{hl}^r + \sum_{i,j \in N} \gamma_{ij}^r y_{ij}^r \quad \forall r \in T \label{eq_ppath_hl:grDefinition} \\
& \sum_{\substack{h \in H \\ \text{if } i \in H}} (x_{ih}^r - x_{hi}^r) + \sum_{i,j \in N}  (y_{ij}^r - y_{ji}^r) = \begin{cases}
 1 & \text{if  } i = or^r \\
-1 & \text{if  } i = de^r \\
0 & \text{otherwise}
\end{cases} \quad \forall i \in N, \forall r \in T \label{eq_ppath_hl:minFlowConstraint} \\
& x_{hl}^r \leq z_{hl} \quad \forall h, l \in H, \forall r \in T \label{eq_ppath_hl:openFacilityOnlyAvailable} \\
& g^r \leq \min\{g(\pi) : \pi \in A^r \cup RP^r \ \land \ Feasible(\pi,{\bf z})\} \quad \forall r \in T_{latent} \label{eq_ppath_hl:ub} \\
& \lambda_\pi^r = {\sc Selected}(\pi,{\bf x^r}) \quad \forall r \in T_{latent}, \pi \in A^r \label{eq_ppath_hl:selected} \\
& z_{hl}, x_{hl}^r, y_{ij}^r \in \{0, 1\} \quad \forall h, l \in H, \forall i, j \in N, \forall r \in T \label{eq_ppath_hl:binaryConstraint_start} \\
& \lambda_{\pi}^r \in \{0, 1\} \quad \forall r \in T_{latent}, \pi \in A^r \label{eq_ppath_hl:binaryConstraint_end}
\end{alignat}
\end{subequations}
\caption{High-Level Formulation of Model {\sc P-Path}.}
\label{fig:ppathhl}
\end{figure}

This section proposes the second path-based formulation, {\sc P-Path}, that
utilizes the sets $RP^r$ and $A^r$ for each trip $r \in T_{latent}$. A
high-level presentation of {\sc P-Path} is presented in Figure
\ref{fig:ppathhl}. The critical novelty in {\sc P-Path} is that
variables $({\bf x^r},{\bf y^r})$ are used for every trip $r \in
T$. In other words, {\sc P-Path} computes the paths for every rider
through constraints~\eqref{eq_ppath_hl:bus_flow_balance},
\eqref{eq_ppath_hl:fixed_arc},~\eqref{eq_ppath_hl:grDefinition},
\eqref{eq_ppath_hl:minFlowConstraint}, and
\eqref{eq_ppath_hl:openFacilityOnlyAvailable}. The role of the sets
$A^r$ and $RP^r$ is to ensure that the selected paths for riders with
choice are optimal solutions to the follower subproblem. This is
achieved in two steps. First, the constraints~\eqref{eq_ppath_hl:ub}
ensures that the selected path for trip $r$ is no worse than all the
feasible paths in $A^r \cup RP^r$. Second, constraints
\eqref{eq_ppath_hl:selected} link decision variables $({\bf x^r},{\bf
  y^r})$ with the paths and their selection variables $\lambda^r_\pi$
used in the objective function. The formulation uses the following
abbreviation
\[
{\sc Selected}(\pi,{\bf x^r}) \equiv \forall (h,l) \in x(\pi): x^r_{hl} = 1 \; \wedge \; \forall (h,l) \notin x(\pi): x^r_{hl} = 0
\]
which indicates that path $\pi$ has been selected for trip
$r$. 

\begin{thm}
Formulation~\eqref{eq:ppathhl} returns a design that is optimal for Formulation~\eqref{eq:bl}.
\end{thm}

\proof{Proof:} The proof shows that every design has the same
objective value in both formulations. Consider a design ${\bf
  \bar{z}}$ and a trip $r \in T_{latent}$. The proof makes a case analysis
for each trip $r \in T_{latent}$. Riders in $T_{core}$ do not raise
any issue as discussed earlier. 

Assume first that trip $r$ has one or more profitable paths in ${\cal
  O}^r_f({\bf \bar{z}})$. Such a path necessarily belongs to $AP^r
\cup RP^r$. By property~\ref{prop:selectA}, if there exists such a
path in $A^r$, the bilevel formulation would select it.  {\sc P-Path}
also selects such a path $\pi$. Indeed, constraints
\eqref{eq_ppath_hl:ub} ensures that $\pi$ has the best objective value
among the feasible paths in $AP^r \cup RP^r$. Moreover, the objective
function drives $\lambda_\pi^r$ to 1 to decrease the objective
function (the third term becomes negative). Furthermore, and importantly,
constraints~\eqref{eq_ppath_hl:selected} ensure that a
single path from $AP^r$ is selected in order to avoid decreasing the
objective with multiple profitable paths.  If only $RP^r$ contains
optimal paths, the follower problem selects a path in $RP^r$ that is
profitable but rejected by riders of trip $r$. Here, constraint~\eqref{eq_ppath_hl:selected}
also ensures that all the feasible paths in $AP^r$ are not selected. 
The objective terms in both the bilevel model and the {\sc P-Path}
formulation are thus the same.

If there is no profitable path for trip $r$ in ${\cal
  O}^r_f({\bf \bar{z}})$, the optimal path
returned by the follower problem must belong to $ANP^r \cup RNP^r$. If
such an optimal path exists in $RNP^r$, by Property
\ref{prop:selectRNP}, the bilevel formulation will choose one of them.
    {\sc P-Path} ensures this automatically since the optimization has
    no incentive to assign the $\lambda$ variables to 1. If no such
    path exists, the {\sc P-Path} model is constrained to select the
    optimal path in $A^r$, which is exactly the choice that is made by
    the bilevel formulation. Again constraints
   ~\eqref{eq_ppath_hl:selected} makes sure that the proper
    $\lambda_\pi^r$ variable is set to 1, all others being zeros.
\Halmos\endproof

\begin{figure}[!t]
\begin{subequations}
\label{eq:ppath}
\begin{alignat}{1}
\min_{z_{hl}} \quad & \sum_{h, l \in H} \beta_{hl} z_{hl} + \sum_{r \in T_{core}} p^r g^r + \sum_{r \in T_{latent}} \sum_{\pi \in A^r} p^r \lambda^r_\pi (g(\pi)  - \varphi) \label{eq_ppath:objective}\\
\text{s.t.} \quad 
& \eqref{eq_ppath_hl:bus_flow_balance}-\eqref{eq_ppath_hl:openFacilityOnlyAvailable} \notag \\
& f_\pi^r \leq z_{hl} \quad  \forall r \in T_{latent}, \pi \in A^r \cup {RP}^r \\
& \sum\limits_{ (h, l) \in x(\pi)} z_{hl} - |x(\pi)| + 1 \leq f_\pi^r \quad \forall r \in T_{latent}, \pi \in A^r \cup {RP}^r \label{eq_ppath:feasible_path_requirement}\\
& f_\pi^r \rightarrow g^r \leq g(\pi) \quad \forall r \in T_{latent}, \pi \in A^r \cup {RP}^r \label{eq_ppath:checkAllFeasiblePaths} \\
& \lambda_\pi^r \leq x_{hl} \quad \forall r \in T_{latent}, \pi \in A^r, (h, l) \in x(\pi) \label{eq_ppath:not_avaialibe_arc_then_cannot_select}\\
& \lambda_\pi^r \leq (1 - x_{hl}) \quad \forall r \in T_{latent}, \pi \in A^r, (h, l) \notin x(\pi) \\
& \lambda_\pi^r \geq \sum\limits_{(h, l) \in x(\pi) } x_{hl}^r - |x(\pi)| - \sum\limits_{(h, l) \notin x(\pi)} x_{hl}^r + 1 \quad \forall r \in T_{latent}, \pi \in A^r \label{eq_ppath:must_select_if_adopting_path} \\
& \eqref{eq_ppath_hl:binaryConstraint_start}-\eqref{eq_ppath_hl:binaryConstraint_end}\notag  \\
& f_{\pi}^r \in \{0, 1\} \quad  \forall r \in T_{latent}, \forall \pi \in A^r \cup {RP}^r 
\end{alignat}
\end{subequations}
\caption{The {\sc P-Path} Formulation.}
\label{fig:ppath}
\end{figure}

\noindent
Constraints~\eqref{eq_ppath_hl:ub} can be linearized easily by
using the same variables $f_\pi^r$ as in model {\sc
  C-Path}. Constraints~\eqref{eq_ppath_hl:selected} can be linearized
with the set of constraints \eqref{eq_ppath:not_avaialibe_arc_then_cannot_select}---\eqref{eq_ppath:must_select_if_adopting_path}.

\noindent
The final formulation for {\sc P-Path} is given in Figure~\ref{fig:ppath},
where the implication constraints can be replaced with their linear big-M
transformations. The number of variables and constraints for {\sc P-Path}
are given by: 
\begin{equation*}
\begin{split}
    \text{\# Variables: }&  {|H|}^2 + ({|H|}^2 + {|N|}^2 + 1)\cdot|T| +  \sum_{r \in T_{latent}}  \bigg( 2|A^r| + |{RP}^r| \bigg) \\
    \text{\# Binary Variables: }&  {|H|}^2 + ({|H|}^2 + {|N|}^2 ) \cdot|T| +  \sum_{r \in T_{latent}}  \bigg( 2|A^r| + |{RP}^r| \bigg) \\
    \text{\# Constraints: }&  |H| + |Z_{fixed}| + (1 + |N| +{|H|}^2) \cdot |T| + \sum_{r \in T_{latent}} \bigg(|A^r| (|H|^2 + 4) + 3|{RP}^r|  \bigg)
\end{split}
\end{equation*}

It is worth pointing out that constraint~\eqref{eq_ppath:not_avaialibe_arc_then_cannot_select} can be replaced by the following constraints to reduce model size:
\begin{equation}
    \lambda_\pi^r \leq f_\pi^r \quad \forall r \in T_{latent}, \pi \in A^r
\end{equation}
This is valid because a path cannot be selected if it is not feasible ($f_\pi^r$ = 0).

\begin{remark}
\label{rmk:lex_for_c_and_p}

To demonstrate the equivalence of the {\sc P-Path} formulation with the original bilevel formulation proposed by \cite{Basciftci2021} with the follower problem
\eqref{eq_bi:lowerLevelProblem_sc} (as presented in Section~\ref{subsect:bilevel}), the lexicographic objective can be integrated into the calculation of $g^r$ and $g(\pi)$ values.
This
can be achieved by replacing $\tau^r_{hl}$ and $\gamma^r_{ij}$ values
with $\hat{\tau}^r_{hl} := M \tau^r_{hl} + t_{hl}' + t_{hl}^{wait}$
and $\hat{\gamma}^r_{ij} := M \gamma^r_{ij} + t_{ij}$, respectively, for a sufficiently large big-$M$ value. 
To adjust the objective function value in~\eqref{eq_ppath:objective}, $\beta_{hl}$ and $\varphi$ values can be further replaced with $\hat{\beta}_{hl} := M \beta_{hl}$ and $\hat{\varphi} := M \varphi$, respectively. 
These modifications apply for the {\sc C-Path} formulation as well to
demonstrate its equivalence with the original bilevel formulation.
\end{remark}

To solve {\sc P-PATH}, commercial solvers can be directly applied. Utilizing a
warm-start approach by giving an initial solution as $z_{hl} = 0$ for
$h,l \in H$ such that $(h,l) \notin Z_{fixed}$ can be beneficial. Additionally,
employing priority branching in the branch and bound tree by
prioritizing network design decisions $\mathbf{z}$ over the remaining
binary variables improve the computational performance.

\begin{remark}
\label{rmk:lazy_const_alg}
For addressing the {\sc P-Path} problem, one might think of applying a solution algorithm, which can incrementally add the necessary paths and constraints for more efficiency. To this end, a lazy constraint-generation algorithm can be designed by first considering a subset of paths from the set $A^r \cup {RP}^r$ for every trip $r \in T_{latent}$, namely ${A}^r_{temp} \cup {RP}^r_{temp}$, and then iteratively incorporating paths to this initial problem. 
In each iteration, for every trip $r \in T_{latent}$, the optimality of the selected path  under the given network design can be evaluated by solving the follower's problem under the given design. If the path is not optimal, then the necessary paths are added to ${A}^r_{temp} \cup {RP}^r_{temp}$ to be considered in the subsequent iterations. 
This step corresponds to evaluating Constraint \eqref{eq_ppath_hl:ub} of the P-PATH formulation and incrementally adding the relevant paths. 
The algorithm stops when the paths selected for every trip $r \in T_{latent}$ correspond to the optimal paths under the transit network design identified by the model.
However, experimental results have shown that such approaches do not offer computational advantages in solving the {\sc P-Path} problem. This situation is due to the bilevel nature of the original problem: feasibility and optimality of the selected paths under the resulting network design need to be ensured. 
A detailed discussion on the algorithm and its computational results can be found in Appendix~\ref{sect:ppath_alg}.
\end{remark}


\section{Preprocessing Techniques} 
\label{sect:computational_enhancements}

This section introduces preprocessing techniques that reduce problem
sizes significantly.

\paragraph{Elimination of Shuttle Trips and Paths:} The ODMTS-DA problem always admits feasible paths for all riders: they include direct shuttle trips and
trips that are only using fixed hub legs in $Z_{fixed}$. Let
$\bar{\Pi}^r$ be the set of these paths and $\bar{g}^r$ be the
minimal cost among them, i.e.,
\[
\bar{g}^r = \min_{\pi \in  \bar{\Pi}^r} g(\pi).
\]
Here, $\bar{g}^r$ provides an upper bound on the weighted cost of the optimal path of trip $r$ under any network design, as addition of hub legs to fixed hub legs in $Z_{fixed}$ can only improve the desired value. 
To this end, the set $L^r$ can be reduced by removing
all shuttle legs whose weighted costs are greater than $\bar{g}^r$, i.e.,
\[
\{ (i,j) \in L^r | \gamma_{ij} > \bar{g}^r \}.
\]
Furthermore, this upper bound can be used to reduce the set $\Pi^r$ for every trip $r \in T_{latent}$ by removing the paths \[
\{\pi \in \Pi^r | g(\pi) >  \bar{g}^r\}.
\]
By applying an analogous relationship, the sets $A^r$ and $RP^r$ can
be significantly reduced as well.  Note that the proposed approach in
preprocessing the sets $\Pi^r$, $A^r$, $RP^r$ is valid under any
choice function (viewed as a black-box). These sets can be further
reduced by leveraging the structure of the choice function, e.g.,
exploiting the trip duration and transfer limits in
\eqref{eq:ChoiceFunctionTimeTransfer}. These reduced sets are referred
as $\Pi^r_{red}$, $A^r_{red}$, $RP^r_{red}$ in the case studies in
Sections~\ref{sect:ypsi_case_study} and~\ref{sect:atl_case_study} and Appendix~\ref{sect:ypsilanti_extra} to
showcase the impact of these preprocessing techniques over different
instances.

\paragraph{Path Assignments:} Some paths are guaranteed to be selected by Models {\sc C-Path} and {\sc P-Path}. 
Let
\[
\underline{g}^r = \min_{\pi \in \Pi^r} g(\pi).
\]
Here, $\underline{g}^r$ provides a lower bound on the weighted cost of the optimal path of trip $r$ under any network design. 
If $\bar{g}^r = \underline{g}^r$ and there exists a path $\bar{\pi} \in \Bar{\Pi}^r$ such that
\[
\argmin\limits_{\pi \in  \bar{\Pi}^r} g(\pi) = \{\bar{\pi} \},
\]
then $\bar{\pi}$ is an optimal assignment for trip $r$.  If
$\bar{\pi} \in A^r$, then this assignment is an adopting path and riders in $r$ are guaranteed to adopt.  This
condition can be generalized when there are multiple optimal paths by
reasoning about the path profitability.  More specifically, if there
exists such a path $\bar{\pi} \in \argmin\limits_{\pi \in \bar{\Pi}^r}
g(\pi) \cap AP^r$, then that rider will adopt the ODMTS under any
network design.  Otherwise, that rider will reject the ODMTS.

\paragraph{Rider Removal:}
A trip $r$ can be removed from Models {\sc C-Path} and {\sc P-Path}
when they are guaranteed i) to adopt a specific path or ii) not to
adopt any path. The first case realizes when there exists an optimal
path assignment from adopting paths as discussed above. This can be
pre-computed and added to the resulting objective of the optimization
models.  The latter case occurs when the optimal path assignment
belongs to set $RP^r \cup RNP^r$. The guaranteed rejection can further
happen when $A^r$ is initially empty or $A^r_{red}$ becomes empty
after preprocessing. After removing these trips from the set $T_{latent}$, the
reduced set is referred to as $T_{latent, red}$.

\paragraph{Removal of Hub Legs:}
It is also possible to remove hub legs that are too far from the
origin and destination of a rider. Let $\underline{g}_{or, h}^r$ be
the weighted cost from origin $or$ to hub $h$ and $\underline{g}_{l,
  de}^r$ be the weighted cost from hub $l$ to destination $de$, where
these costs are computed by identifying the shortest paths between
these locations when all hub legs are available. If
\[
\underline{g}_{or, h}^r + \tau_{hl} > \bar{g}^r \; \vee \; \underline{g}_{l,de}^r + \tau_{hl} > \bar{g}^r,
\]
then variable $x^r_{hl}$ can be removed from Models {\sc C-Path} and
{\sc P-Path}, as its inclusion will not contribute to the objective
function.

Combining these preprocessing techniques, the sizes of the
formulations {\sc C-Path} and {\sc P-Path} can be significantly
reduced, providing computational efficiency.  Under the generalized setting (subproblem \eqref{eq_bi:lowerLevelProblem}), when solving {\sc P-PATH} with linearized constraints, the big-M values introduced to linearize constraint \eqref{eq_ppath:checkAllFeasiblePaths} for each trip $r$ can simply be chosen as $\Bar{g}^r$, given that all paths with greater weighted cost than $\Bar{g}^r$ are already eliminated. 

When enumerating useful paths for {\sc P-PATH}, applying the previously mentioned preprocessing techniques are beneficial. Additionally, leveraging the underlying structure of the choice model $\mathcal{C}^r$ can also be advantageous. In the case studies in Sections~\ref{sect:ypsi_case_study} and~\ref{sect:atl_case_study}, the Algorithm PE-DCM (Path Enumeration Dedicated to Choice Model~\eqref{eq:ChoiceFunctionTimeTransfer}) is applied because the case studies focus on Choice Model~\eqref{eq:ChoiceFunctionTimeTransfer}. This tailored algorithm can efficiently generate $A^r$ and $RP^r \enspace \forall r \in T$ by employing both preprocessing techniques and leveraging the structure of choice model~\eqref{eq:ChoiceFunctionTimeTransfer}. In general, the preprocessing techniques provide a significantly smaller search space for each trip when enumerating its paths for {\sc P-PATH}. Moreover, due to the significance of the transfer tolerance parameter ($l_{ub}^r$) in the choice model, the algorithm has the capability to eliminate all paths containing more than $l_{ub}^r + 1$ arcs. The technical details of this algorithm is outlined in Appendix~\ref{sect:pathEnumeration}.

\section{The Case Study in Ypsilanti}
\label{sect:ypsi_case_study}
This section presents a case study using a realistic dataset from
AAATA\footnote{\url{https://www.theride.org} Last Visited Date:
  January 16, 2024}, the transit agency serving Ypsilanti area and
the broader of Ann Arbor of Michigan, USA. More specifically, Section
\ref{subsect:ypsi_setup} first describes the experimental settings
used in this case study. Section~\ref{subsect:ypsi_exp_result} then
demonstrates the workability and computational advantages of the
path-based formulations by comparing their computational results with 
the exact decomposition algorithm by \citet{Basciftci2021} and the
path-based methods presented in Sections~\ref{sect:paths_method} and
\ref{sect:computational_enhancements}.

\subsection{Experimental Settings}
\label{subsect:ypsi_setup}
\begin{table}
\resizebox{\textwidth}{!}{
\SingleSpacedXI
\begin{tabular}{c c c }
\toprule
Sets and Parameter & Set Size or Parameter Value & Additional Notes \\ \midrule
$N$ & 1267 & Visualized in Figure~\ref{fig:yp_stop_loc}, Appendix~\ref{sect:ypsilanti_extra} \\ 
$H$ & 10 & Visualized in Figure~\ref{fig:yp_stop_loc}, 90 $z_{hl}$ variables need to be determined \\ 
$T$ & 1503 & - \\ 
$T_{core}$ & 1194 (Instances  1 and 3), 937 (Instances  2 and 4) & More details can be found in Table~\ref{table:ypsi_exp_setup} \\
$T_{latent}$ & 309 (Instances  1 and 3), 566 (Instances  2 and 4) & More details can be found in Table~\ref{table:ypsi_exp_setup} \\ 
$Z_{fixed}$ & $\emptyset$ & No backbone lines are considered in Ypsilanti \\
\midrule
\multirow{2}{*}{$n_{hl}$} & \multirow{2}{*}{16} & New buses' headway is assumed to be \SI{900}{s} (\SI{15}{min}) \\ 
& & ODMTS operation lasts for 4 hours\\
$t_{hl}^{wait}$ & \SI{450}{s}  & The expected waiting time computed from buses' headway \\
$b_{dist}$ &  \$3.87 per kilometer & - \\
$w_{dist}$ &  \$1 per kilometer & - \\
$\alpha^r$ & 1.5 and 2.0 for middle and high income households &  Households' income level are presented in Appendix~\ref{sect:ypsilanti_extra} \\
$\phi$ & \$2.5 per rider & Consistent with AAATA's current ticket price \\
\bottomrule

\end{tabular}
}
\caption{The sets and parameter values applied to the Ypsilanti case study.}%
\label{table:ypsi_general_setup}
\end{table}
The experimental settings for this case study are mainly summarized in
Tables~\ref{table:ypsi_general_setup} and
\ref{table:ypsi_exp_setup}. This regular-sized case study is based on
the AAATA transit system that operates over 1,267 ODMTS stops. In
order to design an ODMTS, 10 stops at high density corridors are
selected as ODMTS hubs, and the other stops are only accessible to
shuttles. No specific restriction is applied on new bus arcs;
therefore, the 10 hubs lead to 90 $z_{hl}$ variables. Moreover, all
existing bus lines are assumed to be eliminated, i.e., there is no
backbone lines preserved in this case study. For new bus arcs,
$d_{hl}'$ and $t_{hl}'$ are all assumed to be equal to $d_{hl}$ and
$t_{hl}$. The dataset entails trips between 6 pm and 10 pm, primarily
consisting of commuting trips from work locations to home.

For ridership, the dataset includes 1,503 trips (distinct O-D pairs)
for a total of 5,792 riders. The mode preference of a rider depends on
her income level, that is to say, a rider from a lower income level
has a higher tolerance to travel time. There are 476 low-income, 819
middle-income, and 208 high-income trips with 1,754, 3,316, and 722
riders respectively. The classification of income level are introduced
in Appendix~\ref{sect:ypsilanti_extra}. An $\alpha^r$ value
associated with each choice model (i.e., Equation
\eqref{eq:ChoiceFunctionTime}) is assigned to each class.  Note that
all trips in the low-income class are treated as members of the core
trips set $T_{core}$; hence, no $\alpha^r$ value is required for
them. For the middle-income and the high-income classes, 2.0 and 1.5
are employed as the $\alpha^r$ values, respectively.

To evaluate the performance of the proposed methodology under
different configurations, four instances are generated as described in
Table~\ref{table:ypsi_exp_setup}. In particular, when the ridership is
doubled, the number of riders of each trip is twice as large.  The
core trips percentage parameter for each income level is utilized to
divide the dataset into core trips and latent trips by varied
partitions.

The on-demand shuttle price is set as \$1 per kilometer and each
shuttle is assumed to only serve one passenger.  For buses, the
operating fee is \$3.87 per kilometer and four buses are assumed to
operate between open legs, resulting in an average of 7.5 minutes waiting
time ($t_{hl}^{wait}$). A fixed \$2.5 ticket price that is in line
with the current AAATA system is selected, regardless of the travel
length and multimodality of the trip.  Inconvenience is measured in
seconds, and the inconvenience and cost parameter $\theta$ is set as
0.001.

\begin{table}[!t]
\resizebox{\textwidth}{!}{
\SingleSpacedXI
\begin{tabular}{c c c c c c c c c c c}
       &\multirow{3}{*}{ridership} &\multicolumn{3}{c}{low income core trips} &\multicolumn{3}{c}{medium income core trips} &\multicolumn{3}{c}{high income core trips}\\\cmidrule(lr){3-5} \cmidrule(lr){6-8} \cmidrule(lr){9-11}
       
       &  &\% trips &\# trips & \# riders &\% trips &\# trips &\# riders &\% trips &\# trips & \# riders\\ \midrule 
       
Instance 1 & regular & 100\% & 476 & 1754 & 75\% & 614 & 2842 & 50\% & 104 & 434\\ \midrule

Instance 2 & regular & 100\% & 476 & 1754 & 50\% & 409 & 2262 & 25\% & 52 & 258\\ \midrule

Instance 3 & doubled & 100\%  & 476 & 3508 & 75\% & 614 & 5684 & 50\% & 104 & 868\\ \midrule

Instance 4 & doubled & 100\% & 476 & 3508 & 50\% & 409 & 4524 & 25\% & 52 & 516\\ \bottomrule
\end{tabular}
}
\caption{The experimental setups for the four instances. For doubled ridership, the number of riders for each O-D pair is multiplied by 2. The core trips percentages for each income level are [100\%, 75\%, 50\%] and [100\%, 50\%, 25\%] for low, medium, and high income trips, respectively.}%
\label{table:ypsi_exp_setup}
\end{table}

The proposed approaches are applied to each of the four instances
presented in Table~\ref{table:ypsi_exp_setup}. First, a benchmark run
on the original bilevel formulation (see subproblem
\eqref{eq_bi:lowerLevelProblem_sc}) is carried out with the exact
algorithm proposed in \citet{Basciftci2021}. The exact algorithm is
set to terminate once the gap reaches 0.1\% or the running time
exceeds 6 hours, and the upper bound of its result is reported as the
objective value. Secondly, runs with {\sc P-Path}
Model are applied to those four instances. To compare these
reformulations with the benchmark approach, {\sc
  P-Path} are adjusted to compute the lexicographic optimal appeared
in original subproblem~\eqref{eq_bi:lowerLevelProblem_sc}, as
discussed in Remark~\ref{rmk:lex_for_c_and_p}. On the other hand, an
additional run with {\sc P-Path} is conducted when the generalized
bi-level formulation (see subproblem~\eqref{eq_bi:lowerLevelProblem})
is considered. All models and algorithms were firstly programmed with Python
3.7 then updated to 3.10, and Gurobi 9.5 is selected as the solver. The online repository \cite{guan2024path} contains all the code, along with a sample test case for reference.

\subsection{Computational Results}
\label{subsect:ypsi_exp_result}

\begin{table}[!t]
\resizebox{\textwidth}{!}{
\SingleSpacedXI
\begin{tabular}{c c c c c c c c}
&&&&&& \multicolumn{2}{c}{Path Enum. Time}  (min)  \\ \cmidrule(lr){7-8} 
  {\sc P-Path} & $\sum\limits_{r \in T_{latent, red}} |{RP}_{red}^r |$ &   $\sum\limits_{r \in T_{latent, red}} | A_{red}^r |$   & \# Vars. & \# Binary Vars. & \# Constrs.  & Black-box & Dedicated  \\
\midrule  
  Instance 1   &     200      &    675       &   136,427      &     135,041         &      203,243     & 27.72 &    0.02    \\ \midrule
  Instance 2   &     284       &     1190      &     126,110     &  124,825            & 237,707           &  46.90  &    0.04  \\ \midrule
  Instance 3   &    200    &   675       &   136,427    &     135,041      & 203,243     &  27.61 & 0.02    \\ \midrule
  Instance 4   &    284 &     1190      &     126,110    &  124,825             & 237,707         &  49.02 & 0.04 \\
\bottomrule
\end{tabular}
}
\caption{The Number of Paths, Variables, and Constraints in the {\sc P-Path} Model After Preprocessing.}
\label{table:ypsi_basics_ppath}
\end{table}

This section presents the computational study over the four instances
in Table~\ref{table:ypsi_exp_setup}. Note that instances 1 and 3 were
used in the previous study \citep{Basciftci2021} under the original
bilevel framework, where the authors articulate the benefits of
deploying an ODMTS. This paper only focuses on computational aspects.

\paragraph{Path Enumeration Efficiency}

Table~\ref{table:ypsi_basics_ppath}
presents the model size and reports the computation times for enumerating the paths for {\sc P-PATH}
model. For black-box path enumeration, Algorithm~\mynameref{alg:s1_alg_general_form} (see Appendix~\ref{sect:pathEnumeration})  can be utilized with all preprocessing techniques applied.
The enumeration algorithms may use the choice models as
black-boxes or they can exploit their underlying structures to prune
the search space earlier.  This is especially helpful when using the
{\sc P-PATH} model since it only takes $A_{red}^r$ and ${RP}_{red}^r$
as inputs, and these sets impose strong conditions on the choice
models. In particular, for the choice model
\eqref{eq:ChoiceFunctionTime}, applying k-shortest path algorithms on
graphs in terms of $t^r$ and $g^r$ makes it possible to enumerate
all paths in $A_{red}^r$ and ${RP}_{red}^r$ almost instantly (see Algorithm~\mynameref{alg:s1_alg_choice_form} in Appendix~\ref{sect:pathEnumeration} for detials). Table
\ref{table:ypsi_basics_ppath} describes the significant improvements
in running times when the structure of the choice models is used
during the enumeration processes.  The k-shortest path algorithms used
in this study are by \citet{yen1971finding}.
Corresponding results for the
{\sc C-PATH} model is presented in Appendix~\ref{sect:ypsilanti_extra}.

\paragraph{Computational Efficiency}
\begin{table}[!t]
\resizebox{\textwidth}{!}{
\begin{tabular}{c c c c c c}
& Bilevel Subproblem Type  & Run & Optimality Gap (\%) & Run Time (min) & Objective \\ \midrule
\multirow{3}{*}{Instance 1} & \multirow{2}{*}{Original}   & Benchmark  & 0.65 & 364.10 & 19012.91\\ 
                            &                             & {\sc P-Path} & 0.00 & 2.19 & 19012.91 \\ \cmidrule{2-6}
                            & Generalized & {\sc P-Path} & 0.00 & 2.70 & 19012.91 \\ \midrule 
\multirow{3}{*}{Instance 2} & \multirow{2}{*}{Original}   & Benchmark & 2.48 & 367.52 & 16635.73 \\ 
                            &                             & {\sc P-Path} & 0.00 & 7.80 & 16635.73 \\ \cmidrule{2-6}
                            & Generalized & {\sc P-Path} & 0.00 & 5.04 & 16635.73 \\ \midrule
\multirow{3}{*}{Instance 3} & \multirow{2}{*}{Original}   & Benchmark & 0.99 & 363.22 & 34732.09 \\ 
                            &                             & {\sc P-Path} & 0.00 & 2.62 & 34732.09\\ \cmidrule{2-6}
                            & Generalized & {\sc P-Path} & 0.00  & 1.86 & 34732.09\\ \midrule
\multirow{3}{*}{Instance 4} & \multirow{2}{*}{Original}   & Benchmark & 1.85 & 360.29 & 29962.79  \\ 
                            &                             & {\sc P-Path} & 0.00 & 3.34 & 29962.79 \\ \cmidrule{2-6}
                            & Generalized  & {\sc P-Path} & 0.00 & 2.32 & 29962.79  \\ \bottomrule
\end{tabular}
}
\caption{Computational Efficiency of the {\sc P-Path} Model. Original stands for the original bi-level framework with subproblem~\eqref{eq_bi:lowerLevelProblem_sc}; thus, the lexicographic optimal are considered. Generalized stands for the bi-level framework with subproblem~\eqref{eq_bi:lowerLevelProblem} that is compatible with a black-box choice function.}
\label{table:ypsi_core_results}
\end{table}

Table~\ref{table:ypsi_core_results} summarizes the computation
results, where the benchmark is the approach proposed by
\citet{Basciftci2021}. Under the
original bilevel formulation, {\sc P-PATH} improves the running time
by at least two orders of magnitude. In fact, {\sc P-PATH} finds the
optimal solution in a few minutes, while the benchmark has not reached
the 0.1\% gap within the time limit of 6 hours. These results
demonstrate the workability and computational efficiency of {\sc
  P-PATH}. Under the generalized bilevel formulation, {\sc P-PATH}
finds optimal solutions even faster. Interestingly, the two bilevel
formulations return the same optimal design on this case study.  Note
also that, under generalized bilevel formulation, the discrete
variables $x_{hl}^r$ and $y_{ij}^r$ can be relaxed to be continuous
variables since the subproblem~\eqref{eq_bi:lowerLevelProblem} becomes
totally unimodular when the design decisions are fixed. This further
improves the computational efficiency of the {\sc P-PATH} model.


\paragraph{The Importance of Preprocessing Techniques}

\begin{table}[!t]
\resizebox{\textwidth}{!}{
\SingleSpacedXI
\begin{tabular}{c c c c c c c}
\# Reduction on & $|T_{latent}|$ & 
 $\sum\limits_{r \in T} |L^r|$
 & 
 $\sum\limits_{r \in T_{latent}} |{RP}^r|$ & $\sum\limits_{r \in T_{latent}} |A^r|$  &  $\sum\limits_{r \in (T_{core}) \cup T_{latent, red}}  \sum\limits_{h,l \in H} x_{hl}^r$  \\  
\midrule
  Instances 1 \& 3   &  117  &  20,002 
  &    27       &         289             &        1839    \\ \midrule
  Instances 2 \& 4  &    218  & 20,002  
  &     55      &          498             &       3761    \\
\bottomrule
\end{tabular}
}
\caption{The Impact of Preprocessing on the {\sc P-Path} Model: Reductions of the Set Sizes.}
\label{table:ypsi_preprocessing_impact}
\end{table}

Table~\ref{table:ypsi_preprocessing_impact} demonstrates the
substantial benefits of the preprocessing techniques: it shows the
number of paths that are eliminated for the various sets. The first
observation is that more than 35\% of the latent trips can be
eliminated through path assignments and rider removal. Moreover, the reductions on $A^r$ and $RP^r$ are also helpful, where the size
reduction comes from
two sources: (i) the smaller number of latent trips in $T_{latent, red}$, and
(ii) the reduced number of allowed shuttle arcs in $L_{red}^r$.
Finally, the number of eliminated $x_{hl}$ variables
is non-negligible.  In practice, computing $\bar{g}^r$, reducing
$L^r$, and eliminating $x_{hl}^r$ are performed separately, and the
total computational time is less than 30 seconds. All other set
reductions are carried out within the path enumeration. In summary,
these results highlight the effectiveness of preprocessing and their
importance when solving ODMTS-DA problems.




\section{Large-scale Case Study in Atlanta}
\label{sect:atl_case_study}

This section presents a large-scale case study conducted with travel
data corresponding to a regular workday in Atlanta, Georgia, USA. The
experiment is designed to evaluate the performance of {\sc P-Path} on
large-scale cases. Section~\ref{subsect:atl_exp_setting} explains the
experimental setting, while Section~\ref{subsect:atl_result} presents
the computational results and the benefits of preprocessing
techniques. Appendix~\ref{appendix:atlanta_extra_result} reports
results on the designed ODMTS for those readers interested in the
outcomes of the optimizations. This complements the results on Ypsilanti
detailed in \cite{Basciftci2021}.

\subsection{Experimental Setting}
\label{subsect:atl_exp_setting}

\begin{table}[!t]
\resizebox{\textwidth}{!}{
\SingleSpacedXI
\begin{tabular}{c c c }
Sets and Parameter & Set Size or Parameter Value & Additional Notes \\ \midrule
$N$ & 2426 & Visualized in Figure~\ref{subfig:atl_stops}, Appendix~\ref{sect:atlanta_exp_set}\\ 
$H$ & 58  & Visualized in Figure~\ref{subfig:atl_hubs}, Appendix~\ref{sect:atlanta_exp_set} \\ 
$T$ & 55,871 & - \\ 
$T_{core}$ & 15,478 & -  \\
$T_{latent}$ & 36,283 &  - \\ 
$Z_{fixed}$ & 692 & Derived from four backbone rail lines, see Appendix~\ref{sect:atlanta_exp_set}  \\
\midrule
\multirow{2}{*}{$n_{hl}$} & \multirow{2}{*}{24} & New buses' headway is assumed to be \SI{10}{min} \\ 
& & ODMTS operation lasts for 4 hours\\
$t_{hl}^{wait}$ & \SI{5}{min}  & The expected waiting time computed from buses' headway \\
$b_{time}$ &  \$72.15 per hour & - \\
$w_{dist}$ &  \$0.621 per kilometer  & Equivalent to \$1 per mile \\
$\alpha^r$ & 1.5 &  See choice model~\eqref{eq:ChoiceFunctionTimeTransfer} \\
$l_{ub}^r$  & 2 or 3  & See choice model~\eqref{eq:ChoiceFunctionTimeTransfer} and Table~\ref{table:atlanta_instances} \\ 
$\phi$ & \$2.5 per rider & Consistent with MARTA's current ticket price \\
\bottomrule
\end{tabular}
}
\caption{The sets and values applied to the East Atlanta case study.}%
\label{table:atlanta_general_setup}
\end{table}

This section presents the experimental settings for the case study,
which are mainly summarized in Tables
\ref{table:atlanta_general_setup} and
\ref{table:atlanta_instances}. The Metropolitan Atlanta Rapid Transit
Authority (MARTA) is the major agency that provides transit services
for the Atlanta Metropolis. The test cases all include the MARTA rail
system as the backbone lines, given their importance to the city.  The
operating price for on-demand shuttle and buses are fixed at \$0.621
per kilometer and \$72.15 per hour, respectively. The inconvenience
and cost parameter $\theta$ is fixed at $7.25 / (60 + 7.25)$. Contrary
to Section~\ref{sect:ypsi_case_study}, this case study uses minute as
the unit of time instead of second. The $\theta$ value and the costs
just presented are adopted from an Atlanta-based ODMTS study conducted
by \citet{auad2021resiliency}. There is a unique bus frequency of six
buses per hour, giving an average waiting time ($t_{hl}^{wait}$) of
five minutes. The ODMTS thus operates 24 buses between $h$ and $l$
during the 4 hours time horizon, when the bus arc $(h,l)$ is
selected. The bus arcs must all involve a connection from or to the
rail and are best viewed as rapid bus transit lines expanding the rail
system. The bus travel distance ($d_{hl}'$) and travel time
($t_{hl}'$) are assumed to be equal to $d_{hl}$ and $t_{hl}$, and
these values are obtained using Graph-hopper\footnote{\url{
 https://www.graphhopper.com/} Last Visited Date: January 16,
  2023}. Given that the rail system in the ODMTS is directly adopted
from the existing transit system, its operating costs is omitted
during the design of the ODMTS, since it is a constant.  Thus,
$\beta_{hl} = 0$ if a rail leg connects hubs $h$ and $l$, and its
travel time $t_{hl}'$ is derived from public rail schedules or the
General Transit Feed Specification (GTFS) files. Lastly, shuttles are
allowed to connect hubs in this case study, and each shuttle is
assumed to only serve 1 passenger, i.e., ride-sharing is not
available. A \$2.5 fee is charged for each ODMTS rider, which is the
ticket price of MARTA. Additional information related to the dataset
can be found in Appendix~\ref{sect:atlanta_exp_set}.

\begin{table}
\resizebox{\textwidth}{!}{
\SingleSpacedXI
\begin{tabular}{c c c c c}

           & \# nearby rail hubs & \# $z_{hl}$  & \# undecided $z_{hl}$ &  Transfer Tolerance $l_{ub}^r$ in choice model~\eqref{eq:ChoiceFunctionTimeTransfer} \\ \midrule
Instance 1 & 1           & 732         & 40                     & 2 \\ \midrule

Instance 2 & 1           & 732         & 40                     & 3\\ \midrule

Instance 3 & 2           & 774         & 82                     & 2\\ \midrule

Instance 4 & 2           & 774         & 82                     & 3\\ \midrule

Instance 5 & 3           & 828         & 136                     & 2\\ \midrule

Instance 6 & 3           & 828         & 136                     & 3\\ \bottomrule

\end{tabular}
}
\caption{The six experimental instances considered in the Atlanta case study.}%
\label{table:atlanta_instances}
\end{table}

To test the scalability of {\sc P-Path}, six instances (see Table
\ref{table:atlanta_instances}) were built by controlling two
additional parameters: (i) the number of nearby rail hubs a
bus-only-hub can connect with and (ii) the transfer tolerance
($l_{ub}^r$) of the passengers in the choice model
\eqref{eq:ChoiceFunctionTimeTransfer}. Values 1, 2, and 3 for the
first parameter correspond to 732, 774, and 828 possible hub arcs. The
two values used for $l_{ub}^r$ are 2 and 3. Note that larger transfer tolerance values are generally impractical for transit riders, especially for riders who reside in a metropolitan region. All passengers are assumed
to share the same transfer tolerance. Overall, these two parameters
lead to six instances with different problem sizes and
complexities. In addition, for the adoption factor in the choice
model, a constant value $\alpha^r = 1.5$ is applied to all
passengers. More details related to these two parameters are presented
in Appendix~\ref{appendix:atlanta_extra_result}. 
All models and
algorithms are again implemented and solved using Python 3.7 (later upgraded to 3.10) and Gurobi 9.5 \citep{guan2024path}.


\subsection{Computational Results}
\label{subsect:atl_result}

\begin{table}[!t]
\resizebox{\textwidth}{!}{
\SingleSpacedXI
\begin{tabular}{c c c c c c  c c c}
&&&&&& \multicolumn{1}{c}{Path Enum. Time (min) }   \\ \cmidrule(lr){7-7} 
  {\sc P-Path} & $\sum\limits_{r \in T_{latent, red}} |{RP}_{red}^r |$ &   $\sum\limits_{r \in T_{latent, red}} | A_{red}^r |$   & \# Vars. & \# Binary Vars. & \# Constrs. & Dedicated  \\
\midrule  
Instance 1   &     33      &    7,827      &   14,206,406     &     16,420         &      18,648,588  &   8.85    \\ \midrule
Instance 2   &     33       &   11,115      &      14,311,802    &  22,996            & 21,121,733          & 30.11    \\ \midrule
Instance 3   &    51    &   11,663       &   15,126,181    &     22,152    & 21,620,809   &  9.39  \\ \midrule
Instance 4   &    51 &     17,098      &     15,270,518  &  35,022        & 26,652,670    &  32.16 \\ \midrule
Instance 5   &    78 &     16,042      &     17,228,315   &  32,991       & 27,797,908     &   12.03 \\ \midrule
Instance 6   &    78 &     27,538      &     17,585,987  &  55,979           & 37,357,885     &    44.47 \\ \bottomrule
\end{tabular}
}
\caption{The Number of Paths, Variables, and Constraints in the {\sc P-Path} Model After preprocessing.}
\label{table:atlanta_basics_ppath}
\end{table}

\begin{table}[!t]
\centering
\resizebox{0.8\textwidth}{!}{
\SingleSpacedXI
\begin{tabular}{c c c r r r c}
& & &  \multicolumn{3}{c}{Computational Time (hours)} & \\ \cmidrule(lr){4-6} 
& Bilevel Sub-problem Type  & Run & For 0\% gap & For 1\% gap & For Opt. Sol. & Objective \\ \midrule
Instance 1 & Generalized & {\sc P-PATH}  & 1.94 & 1.93 & 1.93 & 217944.89 \\ \midrule
Instance 2 & Generalized & {\sc P-PATH}  & 2.35 & 2.34 & 2.35 & 218811.60   \\ \midrule
Instance 3 & Generalized & {\sc P-PATH}  & 3.59 & 2.84& 3.04 &  217196.06   \\ \midrule
Instance 4 & Generalized & {\sc P-PATH}  & 6.09 & 4.25 & 5.98 & 218387.18  \\ \midrule
Instance 5 & Generalized & {\sc P-PATH}   & 46.75 & 9.90 & 35.01 & 212917.59 \\ \midrule
Instance 6 & Generalized & {\sc P-PATH}  & 116.48 & 31.52 & 29.98 & 214463.71 \\ \bottomrule
\end{tabular}
}
\caption{Computational Efficiency of the {\sc P-Path} Model. Generalized stands for the bilevel framework with sub-problem~\eqref{eq_bi:lowerLevelProblem} that is compatible with a black-box choice function.}
\label{table:atlanta_compute_result}
\end{table}

\begin{table}[!t]
\centering
\resizebox{0.8\textwidth}{!}{
\SingleSpacedXI
\begin{tabular}{c c c c c c c}

\# Reduction on & $|T_{latent}|$ & 
 $\sum\limits_{r \in T} |L^r|$ &  $\sum\limits_{r \in T_{latent}} |{RP}^r|$ & $\sum\limits_{r \in T_{latent}} |A^r|$  &  $\sum\limits_{r \in T_{core} \cup T_{latent, red}} \sum\limits_{h,l \in H} x_{hl}^r$  \\  
\midrule
Instance 1 & 30,642 & 4,805,467 & 10,175 & 1,117,353 & 14,990,648 \\ \midrule
Instance 2 & 30,501 & 4,805,467 & 10,175 & 2,587,646 & 15,110,231   \\ \midrule
Instance 3 & 30,261 & 4,805,467 & 10,188 & 1,137,236 & 15,931,908   \\ \midrule
Instance 4 & 30,082 & 4,805,467 & 10,188 & 2,621,658 & 16,103,498  \\ \midrule
Instance 5 & 29,235 & 4,805,467 & 10,229 & 1,171,093 & 18,106,680  \\ \midrule
Instance 6 & 27,746 & 4,805,467 & 10,229 & 2,693,458 & 18,521,145   \\ \bottomrule
\end{tabular}
}
\caption{The Impact of preprocessing on {\sc P-Path} Model.}
\label{table:atlatna_preprocessing_ppath}

\end{table}

This section presents the computational results of the six instances
described in Table~\ref{table:atlanta_instances}. They used all the
preprocessing techniques.

\paragraph{Model Size}

The size of the {\sc P-PATH} model over the six instances are reported
in Table~\ref{table:atlanta_basics_ppath}. They show that profitable
rejecting paths (paths in ${RP}_{red}^r$) are rarely observed in
large-scale instances and their numbers remain relatively small as the
problem size grows. This highlights again the benefits of the {\sc
  P-PATH} formulation, since it only relies on a small number of
paths. In addition, greater value for the transfer tolerance parameter
(instances 2, 4, and 6) leads to considerably more complex problems
since the numbers of adopting paths increase significantly.  The path
enumeration leverages the structure of the choice models.  For choice
model~\eqref{eq:ChoiceFunctionTimeTransfer}, bounding the path length
by $l_{ub}^r + 1$ makes it possible to construct $A_{red}^r$ rapidly,
while k-shortest path algorithms can be employed for
${RP}_{red}^r$. Table~\ref{table:atlanta_basics_ppath} shows that the
preprocessing (i.e., reducing path sets and trip sets during path enumeration) takes non-negligible times, yet it is critical to reduce
the sizes of the MIP models. Observe that the numbers of MIP variables
are in the range of 14 to 18 millions and the numbers of constraints
in the range of 18 to 38 millions.

\paragraph{MIP Solving Time}

Table~\ref{table:atlanta_compute_result} reports the computational
efficiency of the {\sc P-PATH} model. The key takeaway is that {\sc
  P-PATH} can solve all instances to optimality, demonstrating its
scalability. The computational times increase significantly as the
problem size grows. Instances 5 and 6 which have more bus arcs are
particularly challenging, but can still be solved optimally.  Table
\ref{table:atlanta_compute_result} also reports the computational time
needed to reach a 1\% optimality gap. The results highlight that the
two largest instances spend the majority of the computation times on
the last 1\% and the optimality proof. {\em These results show that it
  is now in the realm of optimization technology to solve the ODMTS-DA
  problem at the scale of large metropolitan areas.} 
  Moreover, Table~\ref{table:atlanta_compute_result} presents the time for the solver to find the optimal solution as the upper bound during the branch-and-bound procedures. It is evident that the model often reaches the 1\% gaps first.
   {\sc P-PATH} also
provides a benchmark to evaluate the quality of fast heuristics that were proposed in \citet{guan2022}. Additional details regarding the use of {\sc P-PATH} for assessing the quality of heuristic solutions can be found in Appendix~\ref{sect:hrt_vs_ppath}.

\paragraph{Preprocessing Results}

The impact of preprocessing is shown in Table
\ref{table:atlatna_preprocessing_ppath}. After applying path assignments and rider removal, a vast number of trips are
recognized as being ODMTS adoptions or rejections and are removed from $T_{latent}$. Most of the adoptions are
related to local O-Ds in East Atlanta for which direct shuttle paths
can be assigned to the riders. The rejections typically correspond to
long paths that are cost-effective for the agency but inconvenient
for the riders. Table~\ref{table:atlatna_preprocessing_ppath} also
shows that preprocessing eliminates about 50\% of the variables (about
15 to 18 millions of variables). This is not surprising: many hub arcs
are irrelevant for a given rider, since they are far away from the
origin and destination of the trip. 
Both of these elimination methods utilize paths that exclusively involve backbone arcs and shuttles. This preprocessing technique underscores an advantage of ODMTSs. As emphasized by \citet{van2023marta} in a realistic ODMTS pilot, connecting shuttle legs to rails is prominent, and servicing local trips is also a common scenario. These two types of paths can significantly reduce the problem size before solving {\sc P-PATH}. Applying all preprocessing
techniques usually take about 1.5--3 hours for the six instances. The
most time consuming part is variable elimination (see Section~\ref{sect:computational_enhancements}, removal of hub legs) because every
single hub leg needs to be analyzed with the trip set $T_{core} \cup T_{latent, red}$. However, this is still highly beneficial
considering the reduction in the number of variables and the
computational requirements of the MIP model, which is large-scale.


\section{Conclusion}
\label{sect:conclusion} 

This paper considered the ODMTS Design with Adoption problem
(ODMTS-DA) proposed by \citet{Basciftci2021} to capture the latent
demand in on-demand multimodal transit systems. The ODMTS-DA is a
bilevel optimization problem and \citet{Basciftci2021} proposed an
exact combinatorial Benders decomposition to solve it. Unfortunately,
their proposed algorithm only finds high-quality solutions for
medium-sized cities and is not practical for large metropolitan
cities. The main difficulty is in the tension between the design by
the transit agency, which minimizes a combination of cost and
inconvenience, and the choice model of the riders that expresses their
tolerance to, for instance, transit time and the number of transfers.

This paper revisited the ODMTS-DA problem and presented a novel
path-based optimization model, called {\sc P-Path}, aimed at solving
large-scale instances. The key idea behind {\sc P-Path} is to replace
the follower subproblems by enumerating certain types of paths that
capture the essence : adopting paths and profitable paths rejected by
the choice model. The paper showed that, with the help of these two
sets, the bilevel formulation can be replaced by a single-level
optimization which can be formulated as a MIP model. In addition, the
paper presented a number of preprocessing techniques that help in
reducing the sizes of these path sets and hence the number of
variables and constraints in the MIP model. {\sc P-Path} was evaluated
on two comprehensive case studies: the mid-size transit system of the
Ann Arbor -- Ypsilanti region (which was studied in prior work) and
the city of Atlanta.  For the Ann Arbor case study, the experimental
results show that {\sc P-Path} solves the ODMTS-DA instances in a few
minutes: this contrast with the existing approach that cannot solve
the problem optimally after 6 hours. For Atlanta, the results show
that {\sc P-Path} can solve the ODMTA-DA instances optimally: these
instances are obviously out of scope for the existing method. The
resulting MIPs are large-scale (about 17 millions of variables and 37
millions of constraints), but can be solved in a few hours, except for
the large instances which may take longer. These results show the computational benefits of {\sc P-Path}, which provides a
scalable approach to the design of on-demand multimodal transit
systems with latent demand.

As pointed out in Section~\ref{subsect:preliminaries}, the existing ODMTS-DA exhibits several limitations regarding modeling aspects, including the construction of fixed-routes with multiple bus frequencies, balancing shuttle flows, determining the optimal shuttle fleet size, and enabling ride-sharing for shuttles. While these limitations can be addressed in subsequent problems after resolving the ODMTS-DA (i.e., multiple bus frequencies as in \citet{dalmeijer2020transfer}, design with shuttle ride-sharing and fleet-sizing as in \citet{auad2022ridesharing}, and real-time shuttle dispatching with ride-sharing as in \citet{riley2019column}), it remains interesting for future research to explore integrating these elements into the ODMTS-DA framework. This exploration, particularly extending the {\sc P-PATH} approach on shuttle ride-sharing and bus frequencies, could impact on the results of network design.

\section*{Acknowledgments}
This work was partially supported by NSF Leap-HI Grant 1854684 and the Tier 1 University Transportation Center (UTC): Transit - Serving Communities Optimally, Responsively, and Efficiently (T-SCORE) from the Department of Transportation (69A3552047141).

\begin{APPENDICES}

\section{Path Enumeration Algorithms}
\label{sect:pathEnumeration}
\begin{algorithm}[ht]
\SingleSpacedXI
\caption{PE (Path Enumeration)} 
\label{alg:s1_alg_general_form}
    \begin{algorithmic}[1]
    \FORALL{$r \in T_{latent}$}
        \STATE Construct $G^r = (V^r, E^r)$, with arc weight: 
        $\begin{cases}
           \gamma_{ij} \enspace \forall (i,j) \in L^r \text{(shuttles)} \\ 
           \tau_{hl} \enspace \forall (h, l) \in Z \text{(hub-to-hub)}
        \end{cases}$

        
        \STATE With $G^r$, construct set $\Pi^r$ by enumerating all simple paths that starts from node ${or}^r$ and ends at node ${de}^r$, except the invalid paths, i.e., paths that are solely served by two on-demand shuttles.

        
        \STATE Set $A^r = \emptyset$ and ${RP}^r = \emptyset$.

        \FORALL{$\pi \in \Pi^r$}
        
            \STATE Compute a rider's weighted cost with $g(\pi)$. This value is equivalent to the weighed path length obtained from $G^r$.
                    
            \STATE Evaluate $\pi$ with pre-defined choice model ${\cal C}^r (\pi)$.
            \IF{${\cal C}^r (\pi)$ = 1}
                \STATE $A^r = A^r \cup \{ \pi^r \}$. 
            \ELSE
                \IF {$g(\pi) < \varphi$}
                        \STATE ${RP}^r = {RP}^r \cup \{ \pi \}$.
                \ENDIF
            \ENDIF
            
        \ENDFOR
    
    \ENDFOR
    
    \RETURN $\Pi^r$, $A^r$, and ${RP}^r$ for each trip $r \in T_{latent}$.
    \end{algorithmic}
\end{algorithm}

\begin{algorithm}[ht]
\SingleSpacedXI
\caption{PE-DCM (Path Enumeration Dedicated to Choice Model~\eqref{eq:ChoiceFunctionTimeTransfer})}
\label{alg:s1_alg_choice_form}
    \begin{algorithmic}[1]

    \STATE Set $T_{latent,red} = \emptyset$
    \FORALL{$r \in T_{latent} $}

        \STATE Construct $G_{fixed}^r = (N, E_{fixed}^r)$, with edge weight: 
        $\begin{cases}
           \gamma_{ij} \enspace \forall (i, j) \in L^r   \\ 
           \tau_{hl} \enspace \forall (h, l) \in Z_{fixed} 
        \end{cases}$
        
        \STATE Obtain $\bar{g}^r$ by applying shortest-path algorithms on $G_{fixed}^r$ from ${or}^r$ to ${de}^r$.
        \STATE Construct $L_{red}^r$ by excluding unnecessary shuttles arcs with $\bar{g}^r$
        \item[]

        \STATE Construct $G_{full}^r = (N, E_{full}^r)$, with edge weight: 
        $\begin{cases}
           \gamma_{ij} \enspace \forall (i, j) \in L_{red}^r  \\ 
           \tau_{hl} \enspace \forall (h, l) \in Z  
        \end{cases}$

        \STATE Obtain $\underline{g}^r$ by applying shortest-path algorithms on $G_{full}$ from ${or}^r$ to ${de}^r$.
        
        \item[]

        \STATE Reasoning over $\bar{g}^r$, $\underline{g}^r$, and their corresponding paths, decide \textit{rider removal}.
    
        \STATE If \textit{rider removal} detected, iterate to next $r$, else $T_{latent,red} = T_{latent,red} \cup \{r\}$ \\

        \item[]
        \STATE Set $A_{red}^r = \emptyset$ and $P_{red}^r = \emptyset$.

        \STATE Construct $Z_{red}^r$ by excluding unnecessary hub-to-hub arcs with $\underline{g}^r$ (\textit{removal of hub legs})
        
        \item[]
        \STATE Construct $G^r = (V^r, E^r)$, with edge weight: 
        $\begin{cases}
           \gamma_{ij} \enspace \forall (i,j) \in L^r_{red} \\ 
           \tau_{hl} \enspace \forall (h, l) \in Z_{red}^r 
        \end{cases}$
        \STATE Apply $k$-shortest path algorithm on $G^r$ in terms of $g_\pi^r$, iterate through $\pi^r$.
        \WHILE{$g_\pi^r < \varphi$}
            
            \IF{${\cal C} (\pi^r) = 0$ \AND $g_\pi^r \leq \Bar{g}^r$}
            \STATE $RP_{red}^r = RP_{red}^r \cup \{\pi^r\}$.
            \ENDIF
            
            \STATE Iterate to the next path
        \ENDWHILE

        \item[]
        \STATE Construct $G_{time}^r = (V^r, E^r)$, with edge weight: 
        $\begin{cases}
           t_{ij} \enspace \forall (i,j) \in L^r_{red}  \\ 
           t_{hl}' \enspace \forall (h, l) \in Z_{red}^r 
        \end{cases}$
        \STATE Apply $k$-shortest path algorithm on $G_{time}^r$ in terms of $t_\pi^r$, iterate through $\pi^r$.
        \WHILE{$t_\pi^r \leq \alpha^r \cdot t_{cur}^r$}
            
            \IF{${\cal C} (\pi^r) = 1$ \AND $g_\pi^r \leq \Bar{g}^r$}
            \STATE $A_{red}^r = A_{red}^r \cup \{\pi^r\}$.
            \ENDIF
            
            \STATE Iterate to the next path
        \ENDWHILE
        \STATE Apply \textit{rider removal} if $A_{red}^r = \emptyset$


    \ENDFOR
    
    \RETURN $A_{red}^r$, and $RP_{red}^r$ for each trip $r \in T_{latent, red}$
    \end{algorithmic}
\end{algorithm}

This section complements Section~\ref{sect:paths_method} by first presenting a generic path enumeration algorithm, Algorithm~\mynameref{alg:s1_alg_general_form}, and then its enhanced version, which can leverage multiple techniques introduced in Section~\ref{sect:computational_enhancements} to significantly improve its efficiency and produce useful sets such as $T_{latent, red}$, $\Pi_{red}^r$, $A_{red}^r$, and ${RP}_{red}^r$. Specifically, $\bar{g}^r, \underline{g}^r$ can be employed to determine the removal of arcs, riders, paths, and variables, and then $L^r$ and $Z$ could be substituted with $L_{red}^r$ and $Z_{red}^r$, respectively. Note that Algorithm~\mynameref{alg:s1_alg_general_form} is compatible with all types of choice models.

If a specific type of choice model is chosen, e.g., choice model~\eqref{eq:ChoiceFunctionTime}~or~\eqref{eq:ChoiceFunctionTimeTransfer}, alternative forms of Algorithm~\mynameref{alg:s1_alg_general_form} can be designed to speed up the enumeration processes by leveraging the underlying structure of the choice model. To this end,  Algorithm~\mynameref{alg:s1_alg_choice_form} provides a path enumeration algorithm that leverages the benefits of the choice model described in equation~\eqref{eq:ChoiceFunctionTimeTransfer} along with preprocessing techniques implemented before utilizing {\sc P-PATH}. The primary advantage of Algorithm~\textit{PE-DCM} lies in its exclusivity to model {\sc P-PATH}. This aspect eliminates the need to generate the set $\Pi_{red}^r$, a typically time-consuming step.
In Algorithm~\textit{PE-DCM}, an alternative approach to constructing the set $A_{red}^r$ (refer to lines 21-28) involves leveraging the transfer tolerance $l_{ub}^r$ within the choice model. With $G_{time}^r$, this tolerance can be utilized to generate a set of simple paths with a maximum path length of $l_{ub}^r + 1$. Subsequently, the choice model can be employed to determine the adoption status of these paths. This technique can be efficient when the value of $l_{ub}^r$ is reasonable to realistic choices such as 2 or 3.

\section{{\sc C-PATH}: Optimality and Full Formulation}
\label{sect:c_path_additional}

\paragraph{Optimality} By leveraging the high-level formulation of {\sc C-PATH} introduced in Section~\ref{sect:paths_method}, this section first proves the {\sc C-PATH} is able to find the optimality of the ODMTS-DA problem.

\proof{Proof of Theorem \ref{thm:CPathOptimality}:} The proof shows that every design has the same
objective value in both formulations. Consider a design ${\bf
  \bar{z}}$. For each trip $r \in T_{latent}$, formulation
\eqref{eq_cpath_hl} selects a single path $\pi$ that is feasible and
optimal for design ${\bf \bar{z}}$ due to constraints
\eqref{eq_cpath_hl:optimality} and
\eqref{eq_cpath_hl:must_take_1_path}. Hence, both formulation
\eqref{eq_cpath_hl} and the follower problem can select the same path
$\pi$. Moreover the path $\pi$ contributes to the objective
function only if $\pi \in A^r$, in which case ${\cal C}^r(\pi)$
holds. This means that the terms related to each $r \in T_{latent}$ are the same in both
formulations. The result follows from the fact that, for all $r \in T_{core}$,
the constraints and objective terms are the same in both formulations.
\Halmos\endproof

\paragraph{Full Formulation}
The final
formulation for {\sc C-Path} is given in Figure~\ref{fig:c-path},
where the implication constraints can be easily linearized with big-M
transformations.

\begin{figure}[!t]
\begin{subequations}
\label{eq_cpath}
\begin{alignat}{1}
\min_{z_{hl}} \quad & \sum_{h,l \in H} \beta_{hl} z_{hl} + \sum_{r \in T_{core}} p^r g^r + \sum_{r \in T_{latent}} \sum_{\pi \in A^r} p^r \lambda^r_\pi (g(\pi)  - \varphi) 
\label{eq_cpath:obj_func} \\
& \text{(12b)}-\text{(12f)} \notag \\
 &   f_\pi^r \leq z_{hl} \quad \forall r \in T_{latent}, \pi \in \Pi^r, (h,l) \in x(\pi) \\
&   \sum\limits_{(h,l) \in x(\pi)} z_{hl} - |x(\pi)| + 1 \leq f^r_{\pi} \quad \forall r \in T_{latent}, \pi \in \Pi^r \\
& f^r_\pi \rightarrow m^r \leq g(\pi) \quad \forall r \in T_{latent}, \pi \in \Pi^r \\
& \lambda_\pi^r \leq f_\pi^r \quad \forall r \in T_{latent}, \pi \in \Pi^r \\
& \lambda^r_\pi \rightarrow m^r \geq g(\pi) \quad \forall r \in T_{latent}, \pi \in \Pi^r \\
& \text{(12h)}-\text{(12j)} \notag
\end{alignat}
\end{subequations}
\caption{The {\tt C-Path} Formulation.}
\label{fig:c-path}
\end{figure}

The number of variables and constraints for {\sc C-Path}
are given by
\begin{equation*}
\begin{split}
    \text{\# Variables: }& {|H|}^2 + ({|H|}^2 + {|N|}^2 + 1)\cdot |T_{core}| +  \sum_{r \in T_{latent}} \big( 2 |\Pi^r| + 1 \big)\\
    \text{\# Binary Variables: }& {|H|}^2 + ({|H|}^2 + {|N|}^2 )\cdot |T_{core}| +  \sum_{r \in T_{latent}} 2 |\Pi^r|\\
    \text{\# Constraints: }& |H| + |Z_{fixed}| + (1 + |N| +{|H|}^2) \cdot |T_{core}| + \sum_{r \in T_{latent}} \bigg(4|\Pi^r| + 1  + \sum_{\pi \in \Pi^r} |x (\pi)| \bigg)
\end{split}
\end{equation*}

\section{Ypsilanti Case Study: Additional Information}
\label{sect:ypsilanti_extra}
\begin{figure}[htb]
    \centering
    \includegraphics[width=0.75\textwidth]{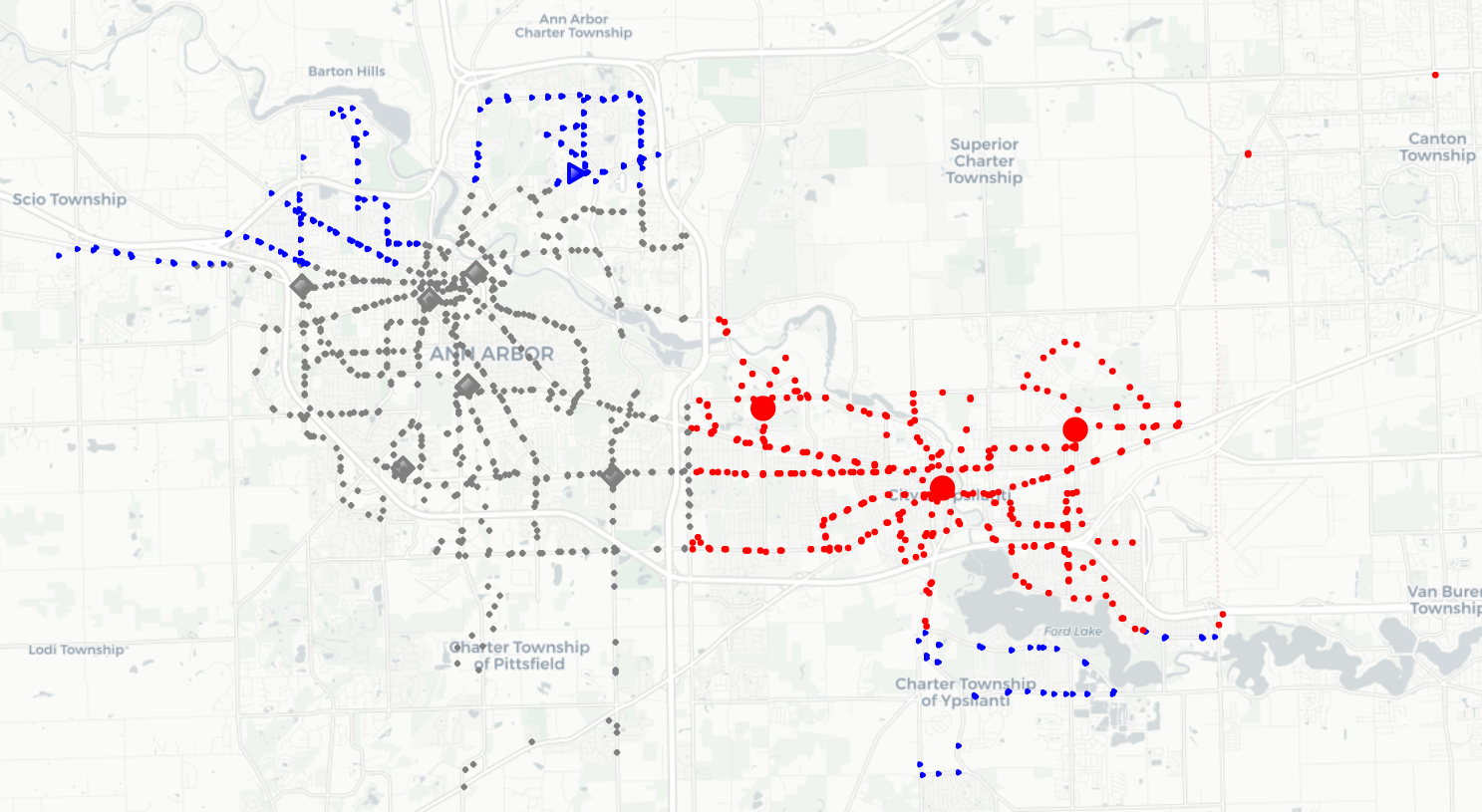}
    \caption{The 1267 bus stops used in this case study. The enlarged 10 stops are the designated hubs. The low-income, middle-income, and high-income stops are represented by red circles, grey squares, and blue triangles respectively. }
    \label{fig:yp_stop_loc}
\end{figure}
\paragraph{Income Level} This section provides additional information for the normal-sized case study. To this end, the whole 1,503 trips were first divided into three different categories: low-income, middle-income, and high-income. The classification of income level is solely decided by a trip's destination because it is approximately a rider's home location. Figure~\ref{fig:yp_stop_loc} visualizes the locations of the stops and their corresponding income levels. All riders belonging to a trip are assumed to have the same income level based on the location of their destination stop.
Following this procedure, there are 476 low-income, 819 middle-income, and 208 high-income trips with 1,754, 3,316, and 722 riders respectively.

\paragraph{Model Sizes Comparison}

\begin{table}[!t]
\resizebox{\textwidth}{!}{
\SingleSpacedXI
\begin{tabular}{c c c c c c c}
&&&&&& Path Enum. Time  (min)  \\ \cmidrule(lr){7-7} 

  {\sc C-Path} & $\sum\limits_{r \in T_{latent, red}} | \Pi_{red}^r |$ &   $\sum\limits_{r \in T_{latent, red}} | A_{red}^r |$   & \# Vars. & \# Binary Vars. & \# Constrs. &  Black-box \\
\midrule  
  Instance 1   &  14,691,440         &    675        &     29,500,425    &   29,499,039           &  155,859,163           &     27.72   \\ \midrule
  Instance 2   &   27,397,108         &       1190    &   54,886,661      &   54,885,376            &    290,162,206        &    46.90    \\ \midrule
  Instance 3   &   14,691,440        &   675         &  29,500,425        &   29,499,039             &      155,859,163       &    27.61   \\ \midrule
  Instance 4   &   27,397,108         &    1190       &   54,886,661      &    54,885,376           &    290,162,206        &     49.02   \\
\bottomrule
\end{tabular}
}
\caption{The Number of Paths, Variables, and Constraints in the {\sc C-Path} Model After Preprocessing.}
\label{table:ypsi_basics_cpath}
\end{table}

In addition to highlighting the importance of using {\sc P-PATH} in practice, the preprocessing techniques and the size of the instances are compared against {\sc C-PATH}.
Table~\ref{table:ypsi_basics_cpath}
summarizes the sizes of model {\sc C-PATH} under the
four instances. Comparing to Table~\ref{table:ypsi_basics_ppath} in Section~\ref{sect:ypsi_case_study}, the substantial difference between the sizes of
{\sc C-PATH} and {\sc P-PATH} models can be observed. In particular, the {\sc P-PATH}
model has two orders of magnitude fewer variables and constraints.
This comes from the sizes of the sets ${RP}_{red}^r$ which are
dramatically smaller than the sets ${\Pi}_{red}^r$. Sets employed by the
{\sc P-PATH} usually consists of a few paths. Optimization results from the {\sc C-PATH} model are not presented, since the MIP solver runs out of memory. However, it is important to recall from Section~\ref{sect:paths_method} that the main role of {\sc C-PATH} is to act as a transitional step for readers between the bilevel formulation and {\sc P-PATH}, rather than being a practical approach for ODMTS-DA. Additionally, the preprocessing techniques remain crucial even if {\sc C-PATH} is applied. The reductions from $\sum\limits_{r \in T_{latent}} |\Pi^r|$ to $\sum\limits_{r \in T_{latent, red}} |\Pi_{red}^r|$ are highly significant. For instances 1 \& 3 and instances 2 \& 4, the amount of reduced paths are 2,900,147,464 and 5,387,002,459, respectively. Note that obtaining all paths in $\bigcup\limits_{r \in T_{latent}} \Pi^r$ requires Algorithm~\mynameref{alg:s1_alg_general_form} to be executed without any preprocessing techniques applied.

\section{Experimental Settings in Atlanta}
\label{sect:atlanta_exp_set}    

\subsection{Data}
\label{subsect:atlanta_data}
This section reveals the details of the data used in the Atlanta case study. The four MARTA rail lines are presented in Figure~\ref{subfig:atl_rails} and all of their 38 stops are reserved as hubs. Therefore, with the method introduced in \citet{dalmeijer2020transfer}, 692 rail arcs can be derived and they form the set $Z_{fixed}$. For $(h, l)$ in $Z_{fixed}$, hub $h$ and $l$ must be taken from the same rail lines (e.g., the blue line), a transfer between rail lines can be represented by concatenating two arcs from $Z_{fixed}$. Moreover, as shown in Figures~\ref{subfig:atl_stops} and~\ref{subfig:atl_hubs}, there are 2,426 ODMTS stops and 58 hubs. The locations of these stops and hubs are selected from a previous study on the ODMTS-DA problem \citep{guan2022}. Because of the existent four rail lines, the 58 hubs can be further classified into two categories with self-explanatory names---20 \textit{bus-only-hubs} and 38 \textit{rail hubs}. Similarly, the demand data used in this section are directly imported from the same study, and these data represent the morning peak demand (6am---10am) of a regular workday. The core demand and latent demand consists of 15,478 and 36,283 riders, and they are visualized in Figure~\ref{subfig:atl_core_trips_heatmap} and Figure~\ref{subfig:atl_additional_trips_heatmap}, respectively. The core data are collected from existing MARTA users and the latent demand data represent drive-alone commuters who reside in East-Atlanta. All details related to data construction can be found in the previous study \citep{guan2022}.
\begin{figure}[!ht]
    \centering
    \begin{subfigure}[b]{0.3\textwidth}
        \includegraphics[width=\textwidth]{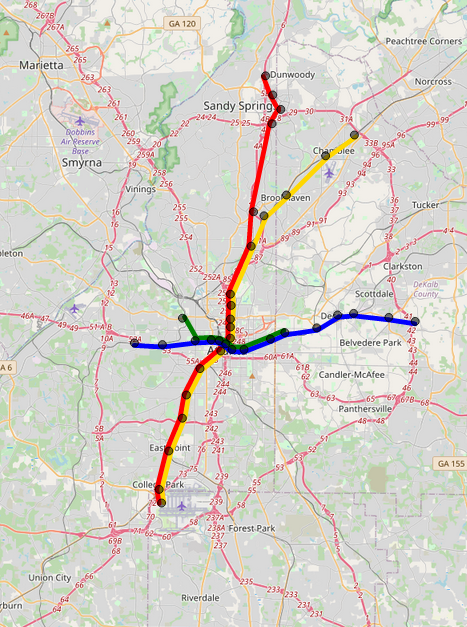}
        \caption{MARTA Rail Lines}
        \label{subfig:atl_rails}
    \end{subfigure}
    \begin{subfigure}[b]{0.3\textwidth}
        \includegraphics[width=\textwidth]{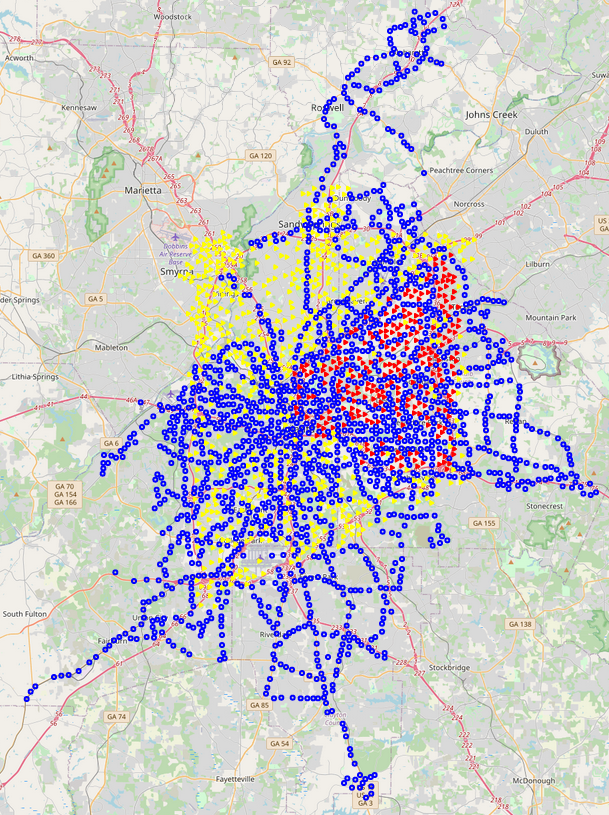}
        \caption{ODMTS stops}
        \label{subfig:atl_stops}
    \end{subfigure}
    \begin{subfigure}[b]{0.3\textwidth}
        \includegraphics[width=\textwidth]{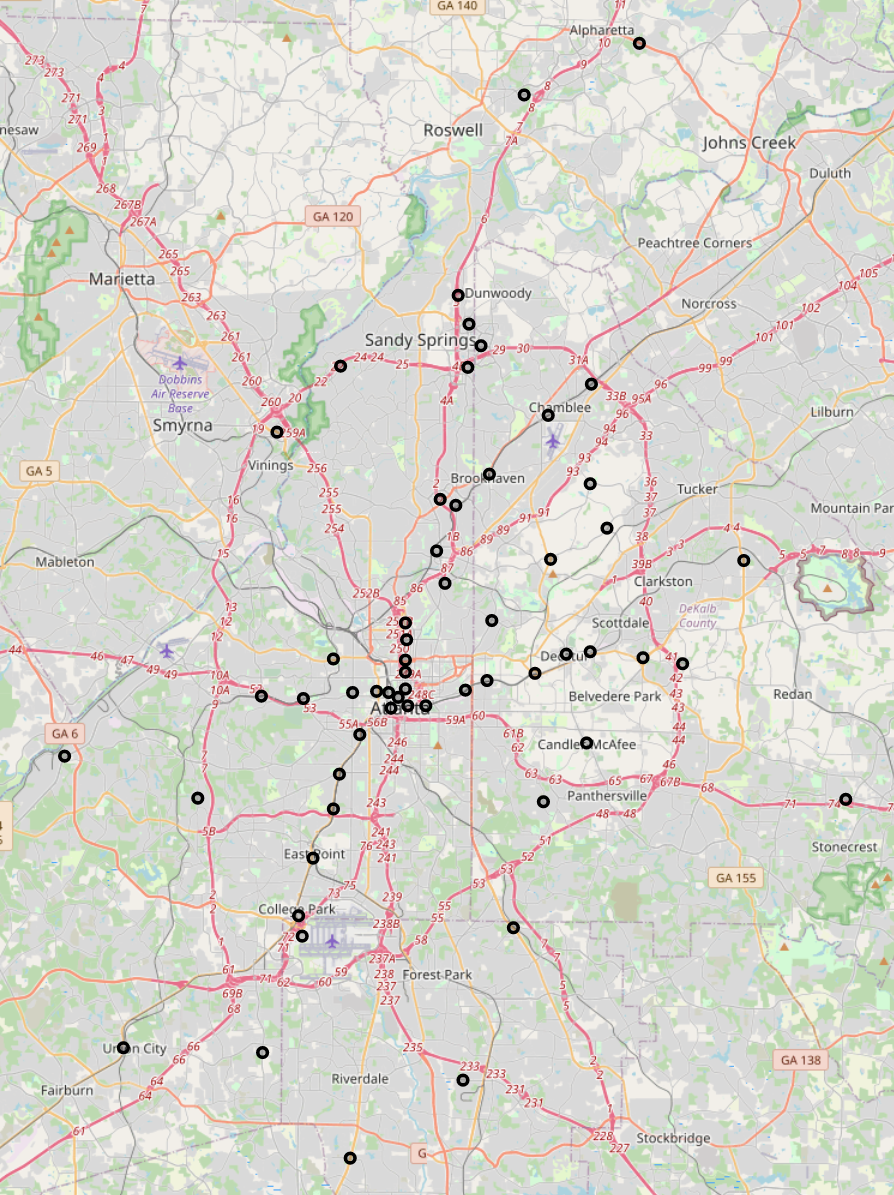}
        \caption{ODMTS Hubs}
        \label{subfig:atl_hubs}
    \end{subfigure}

\caption{The MARTA rail system is shown in (a). The 2,426 ODMTS stops and 58 ODMTS hubs and the MARTA rail system are visualized in (b) and (c), respectively. All 38 MARTA rail stations are reserved as ODMTS hubs.}
\label{fig:atlanta_stop_hub_rail}
\end{figure}
\begin{figure}[!ht]
    \centering 
    \begin{subfigure}[b]{0.45\textwidth}
        \includegraphics[width=\textwidth]{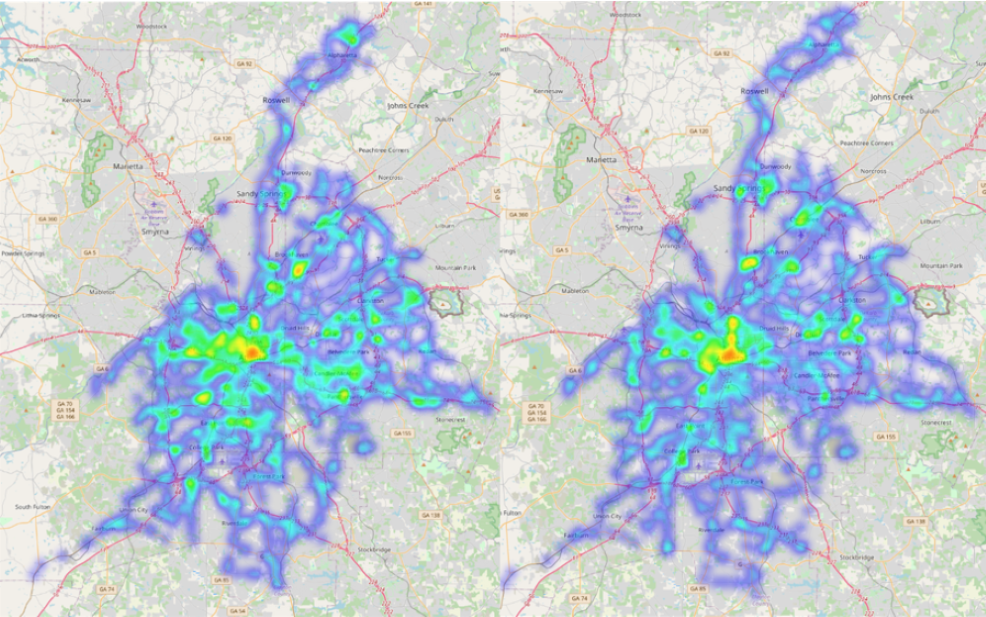}
        \caption{O-Ds in Core Trip Set $T_{core}$.}
        \label{subfig:atl_core_trips_heatmap}
    \end{subfigure}
    \begin{subfigure}[b]{0.45\textwidth}
        \includegraphics[width=\textwidth]{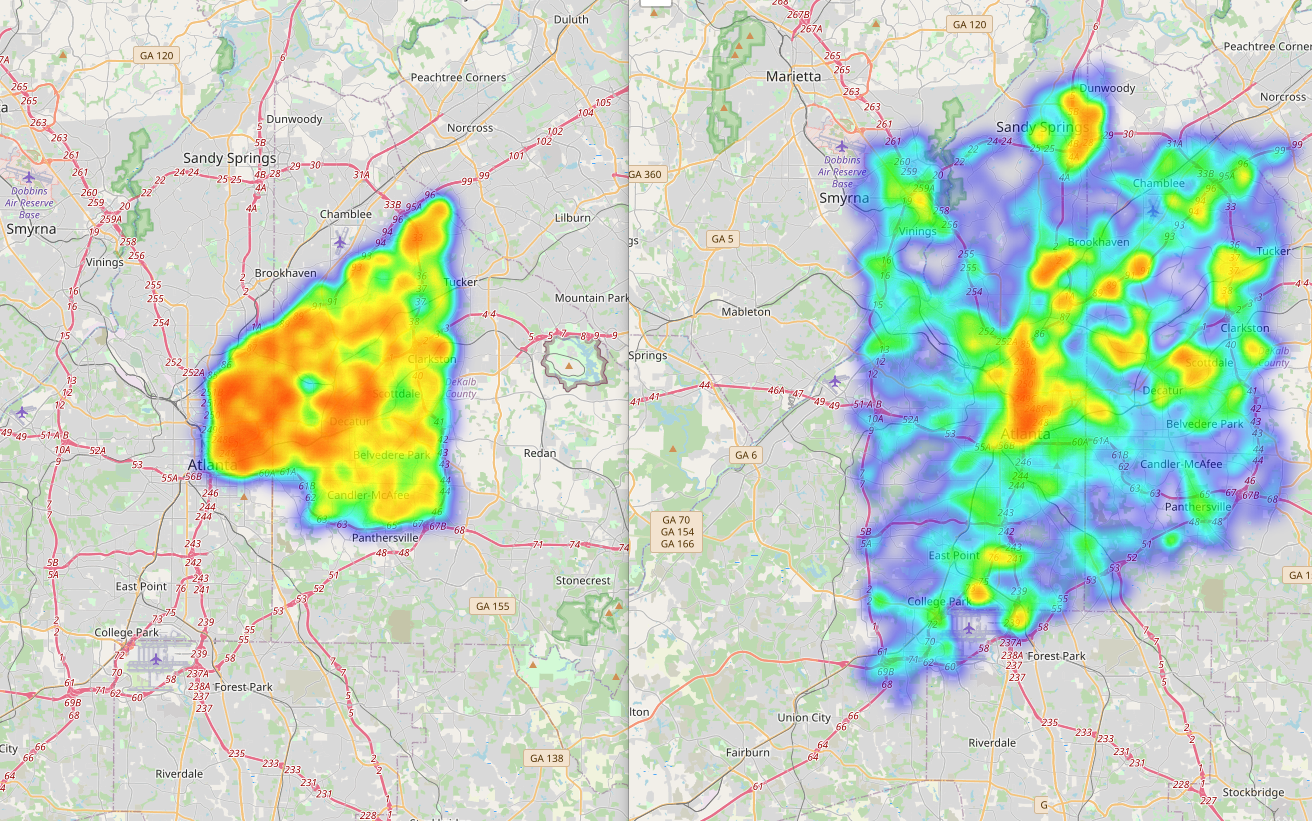}
        \caption{O-Ds in Latent Trip Set $T_{latent}$.}
        \label{subfig:atl_additional_trips_heatmap}
    \end{subfigure}
    
\caption{The heat-maps in (a) and (b) visualize the origin stops and destination stops of the core trip set $T_{core}$ and the latent trip set $T_{latent}$, respectively.}
\label{fig:atlanta_demand}
\end{figure}

\subsection{Parameters}
\label{subsect:atlanta_bus_arc}
This section presents the details of the two parameters used in the Atlanta case study. 
\paragraph{Parameter 1: \# nearby rail hubs }
To select the number of nearby rail hubs for each bus-only-hub, its nearby rail hubs are selected based on the bus travel time $t'_{hl}$. Within these nearby rail hubs, bus connection is allowed such that a $\textit{bus-only-hub} \rightarrow \textit{nearby rail hub 1} \rightarrow \textit{nearby rail hub 2} \rightarrow  \textit{bus-only-hub}$  cycle is possible to be constructed. Additionally, bus arcs are not allowed to overlap rail arcs; in other words, if an arc $(h, l)$ and its symmetric pair $(l, h)$ are elements of $Z_{fixed}$, then they will not be considered as potential bus arcs anymore. Note that the two restrictions above do not contradict each other, it can coexist when the following scenario appears: consider a bus-only-hub, its nearby rail hub 1 is from rail line A and nearby rail hub 2 is from rail line B. Therefore, there is no such arc in $Z_{fixed}$ that already connects nearby rail hub 1 and nearby rail hub 2, then they can potentially be connected with new bus arcs. With above restrictions applied, three values---1, 2, and 3, are chosen, and the sizes of set $Z$ become 732, 774, and 828 in these instances, respectively. 

\paragraph{Parameter 2: transfer tolerance}
In a metropolitan area, a path with more than 3 transfers will hardly be accepted by passengers. Therefore, for the duration and transfer choice model~\eqref{eq:ChoiceFunctionTimeTransfer}, 2 and 3 are tested as transfer tolerance, and all passengers are assumed to share the same values.

\section{Atlanta Case Study: Design Results}
\label{appendix:atlanta_extra_result}

\begin{figure}[!t]
    \centering
    \captionsetup[subfigure]{justification=centering, labelformat=empty}
    \begin{subfigure}[b]{0.3\textwidth}
        \includegraphics[width=\textwidth]{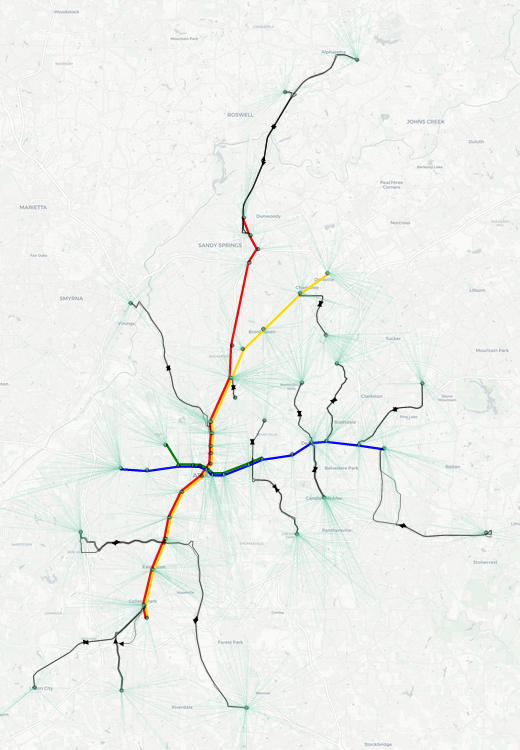}
        \caption{
            Instance 1\\
            1 nearby rail station\\
            transfer tolerance is 2
        }
    \end{subfigure}
    \begin{subfigure}[b]{0.3\textwidth}
        \includegraphics[width=\textwidth]{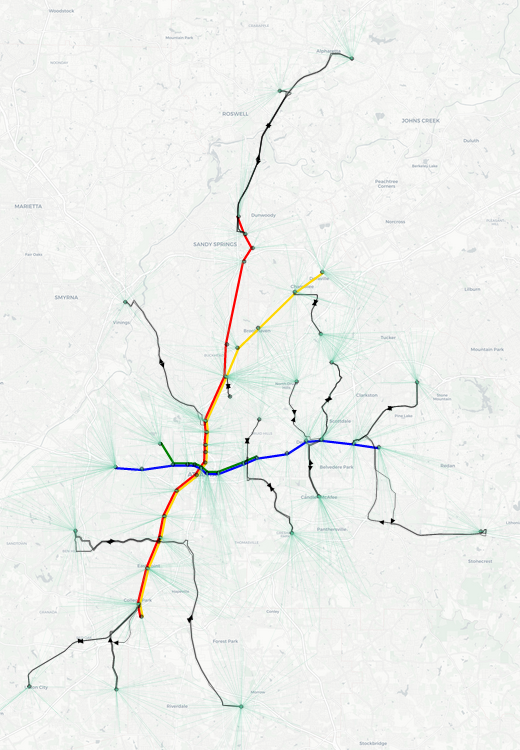}
        \caption{
            Instance 3 \\
            2 nearby rail station\\
            transfer tolerance is 2
        }
    \end{subfigure}
    \begin{subfigure}[b]{0.3\textwidth}
        \includegraphics[width=\textwidth]{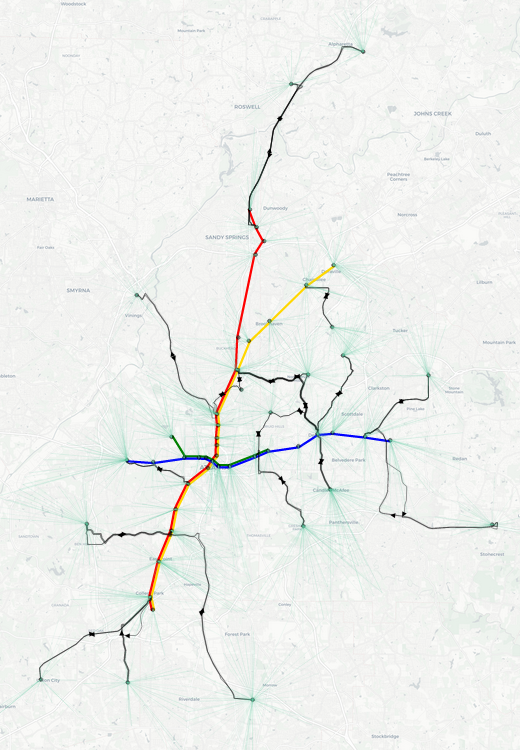}
        \caption{
            Instance 5 \\
            3 nearby rail station\\
            transfer tolerance is 2
        }
    \end{subfigure}
    
    \begin{subfigure}[b]{0.3\textwidth}
        \includegraphics[width=\textwidth]{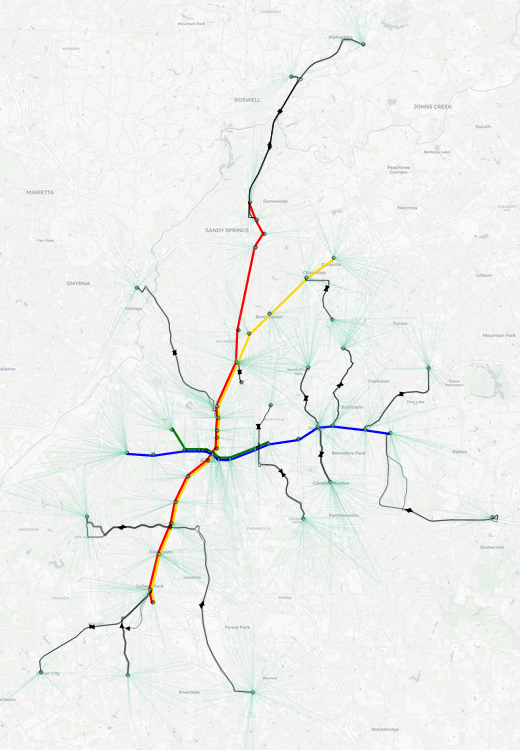}
        \caption{
            Instance 2\\
            1 nearby rail station\\
            transfer tolerance is 3
        }
    \end{subfigure}
    \begin{subfigure}[b]{0.3\textwidth}
        \includegraphics[width=\textwidth]{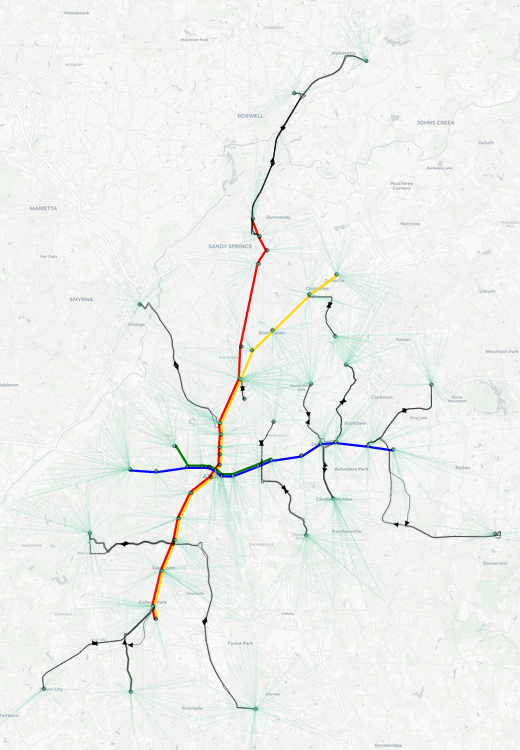}
        \caption{
            Instance 4\\
            2 nearby rail station\\
            transfer tolerance is 3
        }
    \end{subfigure}
    \begin{subfigure}[b]{0.3\textwidth}
        \includegraphics[width=\textwidth]{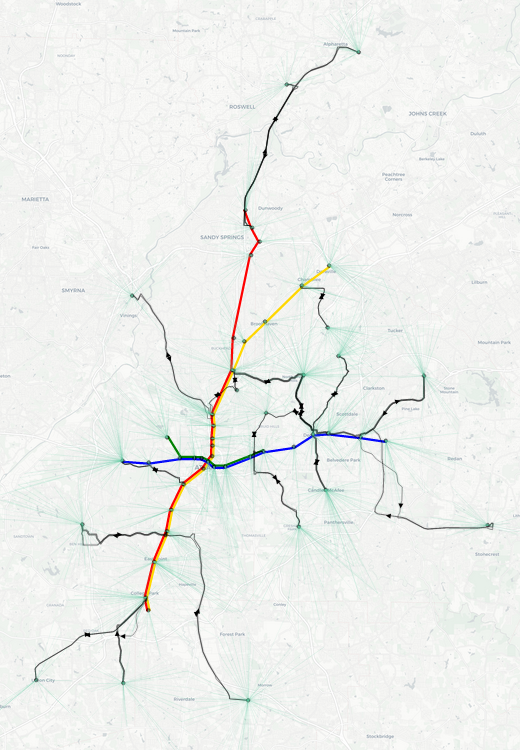}
        \caption{
            Instance 6 \\
            3 nearby rail station \\
            transfer tolerance is 3 
        }
    \end{subfigure}

\caption{The results of all six instances for the large-scale case study are visualized in this figure. The black lines are used to represent the opened new bus arcs decided by the optimization framework. The green lines are used to highlight the non-direct shuttles, i.e., shuttles that connect hubs with trips' origins or destinations.}
\label{fig:atlanta_visual}
\end{figure}

This appendix presents results on the design of ODMTS for Atlanta
using the results of the six instances. It provides the optimal
transit network designs, rider adoptions, and travel duration. The
corresponding optimal designs are visualized in Figure
\ref{fig:atlanta_visual}. It can be seen that designs are similar when
the number of nearby rail stations' is set to 1 or 2. When the number
of nearby rail stations is equal to 3, the designs tend to be more
convoluted in areas with high demand (e.g., East Atlanta) and hub arcs
now link non-connected rail stations.

\begin{table}[!t]%
\centering%
\resizebox{\textwidth}{!}{
\SingleSpacedXI
\begin{tabular}{c c c c c c c c c c c c c c}%

&& \multicolumn{4}{c}{Existing Riders} & \multicolumn{8}{c}{Adopted Riders} \\ %
\cmidrule(lr){3-6} \cmidrule(lr){7-14}

&&\multicolumn{2}{c}{Only Use Shuttle}&\multicolumn{2}{c}{ Use Bus / Rail} &\multicolumn{2}{c}{Adoption}&\multicolumn{2}{c}{Only Use Shuttle}&\multicolumn{2}{c}{ Use Bus / Rail}&\multicolumn{2}{c}{Adoptions with Profits} \\%
\cmidrule(lr){3-4} \cmidrule(lr){5-6} \cmidrule(lr){7-8} \cmidrule(lr){9-10} \cmidrule(lr){11-12} \cmidrule(lr){13-14}

\# nearby rail&tsf tol.&\# riders&\%&\# riders&\%&\# riders&\%&\# riders&\%&\# riders&\%&\# riders&\%\\%
\midrule%
1&2&7292&13.05&48579&86.95&23668&43.11&19989&84.46&3679&15.54&12771&53.96\\%
\midrule%
1&3&7292&13.05&48579&86.95&23852&43.44&19989&83.80&3863&16.20&12794&53.64\\%
\midrule%
2&2&7283&13.04&48588&86.96&23435&42.69&20020&85.43&3415&14.57&12722&54.29\\%
\midrule%
2&3&7283&13.04&48588&86.96&23687&43.14&20020&84.52&3667&15.48&12746&53.81\\%
\midrule%
3&2&7249&12.97&48622&87.03&22382&40.77&18942&84.63&3440&15.37&12736&56.90\\%
\midrule%
3&3&7272&13.02&48599&86.98&22516&41.01&18925&84.05&3591&15.95&12729&56.53\\%
\bottomrule%
\end{tabular}%
}
\caption{Ridership of the reported ODMTS designs.}%
\label{tab:atl_adoption}
\end{table}

Table~\ref{tab:atl_adoption} summarizes the ridership results of the
ODMTS designs. For the existing riders, the majority of them use rails and buses. These can be explained by the fact that existing MARTA users are already heavily relying on the rail system. All ODMTS designs provide at least a 40\% adoption rate, and a higher transfer tolerance only increases the adoption rate
slightly. Unlike the existing riders, a significant percentage only
uses direct shuttle trips among these adopted riders. Many of of these riders are identified as
adoption trips during preprocessing. Although the proportion of
direct shuttle users is high, there is still a large number of riders
using bus or rail services. Recalling that there are 55,871 existing riders in
the system, and these newly adopted bus and rail riders increase the
number of public transit users by at least 6.0\%. Another interesting
observation from Table~\ref{tab:atl_adoption} is at least half of the
newly adopted riders bring profit to the transit agency.

\begin{table}%
\centering%
\resizebox{\textwidth}{!}{
\SingleSpacedXI
\begin{tabular}{c c  c c  c c  c c }%
&&\multicolumn{2}{c}{Riders Adopting ODMTS}&\multicolumn{2}{c}{Existing Riders}&\multicolumn{2}{c}{Riders not adopting ODMTS}\\%
\cmidrule(lr){3-4} \cmidrule(lr){5-6} \cmidrule(lr){7-8}
\# nearby rail&tsf tol.&ODMTS&direct&ODMTS&direct&ODMTS&direct\\%
\midrule%
1&2&7.78&6.81&21.99&13.31&30.66&14.47\\%
\midrule%
1&3&7.96&6.93&21.99&13.31&30.66&14.42\\%
\midrule%
2&2&7.63&6.72&21.99&13.31&30.76&14.48\\%
\midrule%
2&3&7.88&6.88&21.99&13.31&30.76&14.42\\%
\midrule%
3&2&7.61&6.63&21.63&13.31&28.83&14.29\\%
\midrule%
3&3&7.78&6.74&21.79&13.31&29.13&14.25\\%
\bottomrule%
\end{tabular}%
}

\caption{Average travel time of different groups.}%
\label{table:atl_result_travel_time}
\end{table}

Table~\ref{table:atl_result_travel_time} presents the ODMTS travel
time of three different groups: (i) riders from latent trips who adopt the
system, (ii) core riders corresponding to existing transit users, and
(iii) riders from latent trips who continue to drive alone. For each group,
the travel times of two modes---ODMTS and direct driving, are
compared. For the riders who adopt ODMTS, the average ODMTS travel
time is only slightly greater than the average direct driving
time. This can be explained by recalling the direct shuttle percentage
reported in Table~\ref{tab:atl_adoption}. On the other hand, the
non-adopting commuters experience a larger time gap, where the average
duration of the paths suggested by ODMTS are approximately two times
of the average direct trip duration. It is also important to compare
the travel time of the existing system for the core riders. On
average, the riders spend 40 minutes if they only use the existing bus
service, and this value is 26 minutes for existing riders who also use
rails. Therefore, the ODMTS is a more convenient and profitable
choice compared to the existing transit system even the desires of the
latent ridership are considered. 

\begin{table}[!t]%
\centering%
\resizebox{\textwidth}{!}{
\SingleSpacedXI
\begin{tabular}{c c c c c c c}%
\# nearby rail&tsf tol.&Revenue (\$) &New Bus Investment (\$)&Shuttle Cost (\$)&Rail Opt. Cost (\$) &NP/rider (\$)\\%
\midrule%
1&2&198.85k&8.27k&124.52k&24.87k&0.52\\%
\midrule%
1&3&199.31k&8.27k&125.24k&24.87k&0.51\\%
\midrule%
2&2&198.26k&8.33k&123.61k&24.87k&0.52\\%
\midrule%
2&3&198.90k&8.33k&124.63k&24.87k&0.52\\%
\midrule%
3&2&195.63k&10.31k&117.65k&24.87k&0.55\\%
\midrule%
3&3&195.97k&10.16k&118.22k&24.87k&0.54\\%
\bottomrule%
\end{tabular}%
}
\caption{Revenue and cost from the agency perspective. Opt. and NP stand for Operating and Net Profit, respectively.}%
\label{tab:atl_result_cost}
\end{table}

The cost analysis from the perspective of the transit agency is shown
in Table~\ref{tab:atl_result_cost}. Ticket is the only source of
revenue considered in this table, and \textit{ticket revenue} stands
for the revenue collected from the core riders and the newly adopted
riders. The \textit{bus investment} and \textit{shuttle operating
  cost} are directly computed from the equations introduced in Section
\ref{sect:ProblemStatement} without the multiplication of $\theta$ or
$1 - \theta$. This study assumes that the ODMTS is operated by MARTA;
thus, a fixed \textit{rail operating cost} between 6am---10am is also
considered in this table. The rail operating cost is estimated based
on the MARTA rail schedule and the approximated \$240 operating cost
per revenue hour is reported in fiscal year 2018 (shown in Section
4.7.1 in \citet{marta2019report}). After including the 24.87k rail
operating cost, \textit{net profit per rider} values are
computed. Table~\ref{tab:atl_result_cost} shows that the ODMTS is 
currently earning about 50 cents when serving a customer---either a
core rider or a newly-adopted rider. This is a result of the on-demand
shuttles and the rapid bus services efficiently replacing the
traditional buses, especially the buses that operate in regions with
low ridership.

\begin{table}[!t]%
\centering%
\resizebox{\textwidth}{!}{
\SingleSpacedXI
\begin{tabular}{c c c c c c }%
\# nearby rail&tsf tol.&drive alone distance (km)&with ODMTS distance (km)&\% distance reduced&bus travel distance (km)\\%
\midrule%
1&2&576.68k&537.26k&6.83&8.04k\\%
\midrule%
1&3&576.68k&534.51k&7.31&8.04k\\%
\midrule%
2&2&576.68k&539.43k&6.46&8.15k\\%
\midrule%
2&3&576.68k&535.86k&7.08&8.15k\\%
\midrule%
3&2&576.68k&538.54k&6.61&9.79k\\%
\midrule%
3&3&576.68k&536.36k&6.99&9.69k\\%
\bottomrule
\end{tabular}%
}
\caption{Travel Distance of riders with choices. \textit{Drive alone distance} aggregates the travel distance covered by these riders when all of them drive to commute. On the other hand, \textit{with ODMTS distance} sums the adopting riders' travel distance on shuttles and the non-adopting riders' self-driving distances.}%
\label{tab:atl_result_distance}
\end{table}

Lastly, Table~\ref{tab:atl_result_distance} summarizes the ODMTS's
impact on traveled distance by cars. Since core riders are existing
transit users, they are not included in this analysis. There are
substantial reductions in travel distance after a group of
self-driving commuters adopt the ODMTS, which is partly replaced by
the bus or rail services. Table~\ref{tab:atl_result_distance} also
reports the travel distance of the ODMTS buses.

\section{Lazy Constraint-Generation Algorithm for {\sc P-PATH}}
\label{sect:ppath_alg}
\begin{algorithm}[!ht]
\SingleSpacedXI
\caption{Lazy-Constraints Algorithm for {\sc P-PATH}} 
\label{alg:p_path_lazy}
    \begin{algorithmic}[1]

    \REQUIRE The initial set of $A^r_{temp}, RP^r_{temp} \forall r \in T_{latent}$. Default $A^r_{temp} = \emptyset, RP^r_{temp} = RP^r$
    \WHILE{True}
    \STATE Construct {\sc P-PATH-Partial}, replace $A^r$ with $A_{temp}^r$, $RP^r$ with $ RP_{temp}^r \enspace \forall r \in T_{latent}$
    \STATE Solve {\sc P-PATH-Partial} and obtain network design $\mathbf{z}$
    \FOR{$r \in T_{latent}$}
        \STATE Obtain path's weighted cost $g^r_{temp}$ using the solutions of $\mathbf{x}^r$ and $\mathbf{y}^r$

        \FOR{$\pi^r \in (A^r \cup RP^r) \setminus ( A^r_{temp} \cup RP^r_{temp} )$}
            \IF{
                $\pi^r$ is feasible under network design $\mathbf{z}$ 
                and 
                $g^r_{temp} > g_\pi^r$
            }
                \STATE Check the adoption status of $\pi^r$, expand $A_{temp}^r$ or $RP^r_{temp}$ with $\pi^r$
            \ENDIF
        \ENDFOR
    \ENDFOR

    \item[] 
    \IF{$A^r_{temp}$ and $RP_{temp}^r$ do not grow $\forall r \in T_{latent}$}
        \STATE break while loop
    \ENDIF
    
    \ENDWHILE
    \end{algorithmic}
\end{algorithm}

This section completes the discussion in Remark~\ref{rmk:lazy_const_alg} by presenting a Lazy-Constraints Algorithm and its computational results. Due to the bilevel nature of the ODMTS-DA problem, the algorithm needs to ensure optimality of the selected paths under the resulting network design, while considering all feasible paths from $A^r \cup {RP}^r$ and ensuring that the network design identified gives the optimal objective function value for the transit agency's objective. Otherwise, additional paths need to be incorporated into this problem and the problem needs to be resolved with this extended subset of paths, resulting in an iterative solution approach. 
\\
\indent In particular, Algorithm~\ref{alg:p_path_lazy} describes the procedures that are employed in this study. Instead of directly solving {\sc P-PATH}, {\sc P-PATH-Partial} is constructed with $A^r_{temp}$ and $RP^r_{temp} \enspace \forall r \in T_{latent}$. At steps 6---9, the algorithms checks if the founded paths violate the Constrains~\eqref{eq_ppath_hl:ub} in the original {\sc P-PATH} formulation. Alternatively, one can choose other methods to expand paths for $A^r_{temp}$ and $RP^r_{temp}$, such as only selecting the path with the lowest weighted cost that violates Constraints~\eqref{eq_ppath_hl:ub}.
In addition, note that the algorithm starts with $RP^r_{temp} = RP^r$ by default because $RP^r$ is usually a small set in practice, as indicated by the case studies.
\\
\indent In Table~\ref{table:p_path_alg}, the computational results over six instances show that the run time of directly solving P-PATH outperforms the Lazy-Constraints algorithm. These two Atlanta instances are the two smallest instances from Section~\ref{sect:atl_case_study}.
Indeed, depending on the initial subset of paths considered, this algorithm may require either many iterations for convergence with longer overall run times or fewer iterations with longer computational times per each iteration.

\begin{table}[ht]
\centering%
\begin{tabular}{l r r r}%
 Instance & Run Time Directly Solve {\sc P-PATH} & Run Time with Algorithm~\ref{alg:p_path_lazy} & \# Iterations 
 \\
 \midrule
 Ypsilanti 1 & 2.70 minutes & 4.42  minutes & 3 
 \\ \midrule
 Ypsilanti 2 & 5.04 minutes & 14.29 minutes & 4
 \\ \midrule
 Ypsilanti 3 & 1.86 minutes & 4.70 minutes & 4 
 \\ \midrule
 Ypsilanti 4 & 2.32 minutes & 5.79 minutes & 5 
 \\ \midrule
 Atlanta 1 & 1.94 Hours & 19.47 hours & 14
 \\ \midrule
 Atlanta 2 & 2.35 Hours & 23.66 hours  & 13
 \\ 
\bottomrule
\end{tabular}%
\caption{The run time of Lazy-Constraints algorithm compared to the directly solve P-PATH. }
\label{table:p_path_alg}
\end{table}
\section{{\sc P-PATH} as Benchmark for ODMTS-DA Heuristics}
\label{sect:hrt_vs_ppath}

\begin{table}%
\centering%
\begin{tabular}{l c c c r r r }%
method&step size&trip rules&\# total itr.&run. time (hour)&design obj&\% opt. gap\\%
\midrule%
{\sc P-PATH} & - & - & - & 1.94 & 217944.89 & 0.00 \\ \midrule
Alg. $\rho$-GRAD &3000&{-}&7&0.06&217944.89&0.00\\%
\midrule%
Alg. $\eta$-GRRE &3000&{-}&6&0.07&223988.10 & 2.77\\%
\midrule%
Alg. $\rho$-GAGR &2000&-&58&0.49&219977.73& 0.93\\%
\midrule%
Alg. arc-s1&{-}&(a)& 17& 0.51& 217944.89 & 0.00\\%
\midrule%
Alg. arc-s1 &{-}&(d)&11&0.24&242375.77 & 11.21\\%
\midrule%
Alg. arc-s2 &{-}&(c), (a)&18&0.44&217944.89 & 0.00\\%
\midrule%
Alg. arc-s2&{-}&(d), (a)&18&0.36&217944.89 & 0.00\\%
\bottomrule
\end{tabular}%
\caption{A comparison between {\sc P-PATH} and the heuristic algorithms on the Atlanta instance 1. {\sc P-PATH} is solved to optimal while the heuristic can provide rapid approximation on the optimal solution. "-" stands for inapplicable elements.}
\label{table:p_path_vs_heuristic_i1}
\end{table}

\begin{table}%
\centering%
\begin{tabular}{l c c c r r r }%
method&step size&trip rules&\# total itr.&run. time (hour)&design obj&\% opt. gap\\%
\midrule%
{\sc P-PATH} & - & - & - & 116.48 & 214463.71 & 0.00 \\ \midrule
Alg. $\rho$-GRAD &3000&{-}&8&0.07&214906.96&0.21\\%
\midrule%
Alg. $\eta$-GRRE &1000&{-}&17&0.15&221265.96&3.17\\%
\midrule%
Alg. $\rho$-GAGR &2000&{-}&69&0.66&214906.96&0.21\\%
\midrule%
Alg. arc-s1&{-}&(a)&20&0.72&214846.58&0.18\\%
\midrule%
Alg. arc-s1 &{-}&(d)&13&0.31&242229.26&12.95\\%
\midrule%
Alg. arc-s2 &{-}&(c), (a)&22&0.60&214846.58&0.18\\%
\midrule%
Alg. arc-s2&{-}&(d), (a)&22&0.50&214906.96&0.21\\%
\bottomrule
\end{tabular}%
\caption{A comparison between {\sc P-PATH} and the heuristic algorithms on the Atlanta instance 6. {\sc P-PATH} is solved to optimal while the heuristic can provide rapid approximation on the optimal solution. "-" stands for inapplicable elements.}
\label{table:p_path_vs_heuristic_i6}
\end{table}

This section contrasts the computational outcomes of {\sc P-PATH} with those five {\sc ODMTS-DA} heuristic algorithms introduced in \citet{guan2022}. The heuristics outlined in \cite{guan2022} primarily involve two iterative components: (1) solving an ODMTS design problem without the consideration of latent demand and (2) adjusting the demand by assessing the designed ODMTS using a mode choice model. The primary goal of these heuristics is to efficiently approximate the optimal solution for large-scale ODMTS-DA problems.

Tables~\ref{table:p_path_vs_heuristic_i1} and \ref{table:p_path_vs_heuristic_i6} provide numerical comparisons for instances 1 and 6 of the Atlanta Case study. The step size and trip rules play a role as specific input parameters that govern the behaviors of the algorithms. Comprehensive insights into these parameters can be found in \cite{guan2022}. Notably, for both instances, the primary observation is, the heuristics yield high-quality ODMTS-DA solutions within remarkably short computational time. Particularly in Instance 1, some heuristics can discover optimal solutions, although this optimality requires confirmation from {\sc P-PATH}.  Besides the inherent algorithmic design, the rapid termination speed can also be attributed to the nature of the optimization problems included in the heuristics at each iteration. To expedite the algorithms, multi-threading can be employed to run the sub-problems at each iteration for the heuristic algorithms, the results in both tables were generated utilizing 12 threads. Furthermore, as {\sc P-PATH} has demonstrated its capability to attain optimal solutions for extensive ODMTS-DA instances, it is now regarded as a benchmark for solving the ODMTS-DA problem. Thus, the availability of this exact approach to evaluate the quality of the solutions approximated by the heuristic algorithms (which are more practical approaches for large-scale instances) is necessary.

\begin{table}[!ht]
\centering%
\small
\begin{tabular}{l r r}%
 Instance & Run Time (hours) Directly Solve {\sc P-PATH} & Run Time (hours) with Warm-Starts \\
 \midrule
 Atlanta 1 & 1.94 & 1.17 \\ \midrule
 Atlanta 2 & 2.35 & 1.70 \\ \midrule
 Atlanta 3 & 3.59 & 2.69  \\ \midrule
 Atlanta 4 & 6.09 & 6.15 \\ \midrule
 Atlanta 5 & 46.75 & 62.90 \\ \midrule
 Atlanta 6 & 116.48 & 87.59 \\
\bottomrule
\end{tabular}%
\caption{The run time of solving {\sc P-PATH} to optimal when incorporating heuristic solutions as warm-starts. The heuristic time is included in the warm-start column. }
\label{table:atlanta_warm_start}
\end{table}

Additionally, it is also worth exploring the potential of using the solutions of the heuristic algorithms as warm-starts for large-scale {\sc P-PATH} models, i.e., the Atlanta instances in the Section~\ref{sect:atl_case_study}. In order to select the solutions to be utilized in the warm-start, each of the proposed heuristics in \cite{guan2022} was first conducted for a duration of 30 minutes. The best solution, which corresponds to the lowest ODMTS-DA objective value, is chosen as the warm-start for {\sc P-PATH}. However, as reported in Table~\ref{table:atlanta_warm_start}, the findings suggest that utilizing these warm-starts in commercial solvers such as Gurobi does not necessarily ensure a faster way in reaching to the optimal solution. While employing warm-starts during branch-and-bound procedures can establish a favorable upper bound, it is the limitation of the lower bound that contributes to the significant duration of the entire solving processes. 
\end{APPENDICES}

\theendnotes





\bibliographystyle{informs2014} 
\bibliography{ref} 

\end{document}